\documentclass[11pt, oneside]{article}
\usepackage[a4paper, portrait, margin=1.1811in]{geometry}
\usepackage[english]{babel}
\usepackage[utf8]{inputenc}
\usepackage[T1]{fontenc}
\usepackage{lmodern}
\usepackage{helvet}
\usepackage{etoolbox}
\usepackage{graphicx}
\usepackage{titlesec}
\usepackage{caption}
\usepackage{booktabs}
\usepackage{geometry}   
\usepackage{amstext,amsfonts,a4}
\usepackage{amsmath,amsthm}
\usepackage{amssymb}
\usepackage[colorlinks, citecolor=cyan]{hyperref}
\usepackage{caption}
\graphicspath{ {./images/} }
\usepackage{scrextend}
\usepackage{fancyhdr}
\usepackage{dsfont}
\usepackage{algorithm}
\usepackage{algorithmic}
\usepackage{multirow}
\usepackage{graphicx}
\usepackage{parskip}
\usepackage[normalem]{ulem}
\usepackage{url}
\usepackage{tabularx}

\usepackage{babel}
\usepackage[usenames,dvipsnames]{color}
\usepackage{listings}
\usepackage{listingsutf8}
\usepackage{color}
\usepackage{pgf,tikz}
\usepackage{afterpage}
\usepackage{comment}

\numberwithin{equation}{section}

\makeatletter
\patchcmd{\@maketitle}{\LARGE \@title}{\fontsize{16}{19.2}\selectfont\@title}{}{}
\makeatother

\newcommand{\Wperp}{W_\perp}

\definecolor{darkgreen}{rgb}{0.0, 0.5, 0.0}
\definecolor{violet}{rgb}{0.580,0.,0.827}

\usepackage{authblk}

\newtheoremstyle{exampstyle}
  {3pt} 
  {3pt} 
  {} 
  {} 
  {\bfseries} 
  {.} 
  {.5em} 
  {} 

\theoremstyle{exampstyle}
\newtheorem{Thm}{Theorem}[section]
\newtheorem{Cor}{Corollary}[section]
\newtheorem{Def}{Definition}[section]
\newtheorem{Lemma}{Lemma}[section]
\newtheorem{example}{Example}[section]
\newtheorem{remark}{Remark}[section]

\DeclareMathOperator{\dist}{dist}

\newsavebox\affbox
\author[1]{\textbf{J. Bisch}}
\author[2]{\textbf{A. Hannukainen}}
\affil[1,2]{Aalto University
Department of Mathematics and Systems Analysis
P.O. Box 11100
FI-00076 Aalto
Finland.
}

\titlespacing\section{0pt}{12pt plus 4pt minus 2pt}{8pt plus 2pt minus 2pt}
\titlespacing\subsection{12pt}{12pt plus 4pt minus 2pt}{8pt plus 2pt minus 2pt}
\titlespacing\subsubsection{12pt}{12pt plus 4pt minus 2pt}{8pt plus 2pt minus 2pt}

\titleformat{\author}{\normalfont\fontsize{10}{15}\bfseries}{\thesection}{1em}{}

\title{\textbf{\huge A Ritz method for solution of  parametric generalized EVPs}}
\date{}

\renewcommand{\thepage}{\arabic{page}}
\begin{document}
    
\maketitle
	
\noindent\rule{15cm}{0.5pt}
\begin{abstract}
This work deals with approximate solution of generalized eigenvalue problem with coefficient matrix that is an affine function of d-parameters. The coefficient matrix is assumed to be symmetric positive definite and spectrally equivalent to an average matrix for all parameters in a given set.    

We develop a Ritz method for rapidly approximating the eigenvalues on the spectral interval of interest $(0,\Lambda)$ for given parameter value. The Ritz subspace is the same for all parameters and it is designed based on the observation that any eigenvector can be split into two components. The first component belongs to a subspace spanned by some eigenvectors of the average matrix. The second component is defined by a correction operator that is a d + 1 dimensional analytic function. We use this structure and build our Ritz subspace from eigenvectors of the average matrix and samples of the correction operator. The samples are evaluated at interpolation points related to a sparse polynomial interpolation method. We show that the resulting Ritz subspace can approximate eigenvectors of the original problem related to the spectral interval of interest with the same accuracy as the sparse polynomial interpolation approximates the correction operator. Bound for Ritz eigenvalue error follows from this and known results. Theoretical results are illustrated by numerical examples. The advantage of our approach is that the analysis treats multiple eigenvalues and eigenvalue crossings that typically have posed technical challenges in similar works. 
  
		\let\thefootnote\relax\footnotetext{
			\small $^{*}$ \textit{
				\textit{E-mail addresses: \color{cyan}joanna.bisch@aalto.fi, hannukainen.antti@aalto.fi}}\\}
		
  \vspace{0.3cm}
		\textbf{\textit{Keywords}}: \textit{subspace method, model order reduction, reduced basis methods, eigenvalue problem}
	\end{abstract}
\noindent\rule{15cm}{0.4pt}

\section{Introduction}


Let $\mathbb{S}^{n\times n}\subset \mathbb{R}^{n\times n}$ be the space of symmetric $n\times n$-matrices, $\mathbb{S}^{n\times n}_{++} \subset  \mathbb{S}^{n\times n} $ the subspace of symmetric positive definite matrices, $\Lambda > 0$, $S \subset \mathbb{R}^d$ a bounded set of admissible parameter values and $M\in \mathbb{S}^{n\times n}_{++}$. In addition, let $A:\mathbb{R}^d\rightarrow \mathbb{R}^{n\times n}$ satisfy $A(\sigma)=A_0+\sum_{m=1}^d \sigma_m A_m$ for all $\sigma\in \mathbb{R}^{d}$ with $\{A_m\}_{m=0}^d\subset \mathbb{S}^{n\times n}$ and $A(S)\subset\mathbb{S}^{n\times n}_{++}$. In this work, we are concerned about solving the  generalized eigenvalue problem (GEVP): for all given $\sigma \in S$ find $(\lambda(\sigma),x(\sigma)) \in (0,\Lambda) \times \mathbb{R}^n$ satisfying
\begin{equation}
\label{eq:evp}
    A(\sigma) x(\sigma) = \lambda(\sigma) M x(\sigma).
\end{equation}
That is, we are interested in finding the eigensolutions $(\lambda(\sigma),x(\sigma))$ with $\lambda(\sigma)$ in the spectral interval of interest $(0,\Lambda)$ and $x(\sigma)\in\mathbb{R}^n$. 

Multiparametric GEVP \eqref{eq:evp} arises in multiple contexts, e.g., FE-discretization of elliptic PDEs with spatially varying coefficient functions \cite{Rodriguez17,Gardini19} or mechanical vibration problems \cite{Verhoosel06,Hakula15,Hakula19}.  The eigenvalues of GEVPs play an important role in a very large number of applications, for example, in civil engineering or aeronautics \cite[Section 3.5.2]{Chatelin12}. 
In these applications, one is typically interested in finding few of the  smallest eigenvalues and their corresponding eigenvectors for many $\sigma \in S$.

One could be tempted to solve problem \eqref{eq:evp} for each parameter $\sigma\in S$ by applying standard numerical solvers such as Lanczos method \cite{Jones89}, Arnoldi or Jacobi-Davidson \cite{Joost08,Sleij96}, or even by applying a Rayleigh-Ritz method \cite{Sakurai07}. Matrices involved in typical applications leading to \eqref{eq:evp} are often sparse and very large, thus solving even a single GEVP is already expensive. Thus, using traditional approach many times for several values of $\sigma \in S$ is not recommended as the procedure will become prohibitively expensive. To overcome these challenges, various numerical methods for solving multiparametric eigenvalue problems have been developed. Methods based on a Smolyak-Galerkin approximation rely on a sparse grid collocation operator introduced in \cite{Bieri09}. Such methods have been extensively analyzed by R. Andreev and Ch. Schwab \cite{AndSchw11} and have been shown to have higher accuracy compared to traditional Monte-Carlo methods \cite{Bieri09}. 
Low-rank methods have also been exploited to solve stochastic eigenvalue problems \cite{Su19,Benner19} and parameterized linear systems \cite{Kressner11} by constructing low-rank approximation with help of a GMRES-like method \cite{Ballani13} or Richardson-like method \cite{Schwab11}.
Treating crossing or multiple eigenvalues represent a challenge when analyzing the aforementioned methods. 
To ensure convergence rate, authors either do not address such case \cite{HakulaLaaksonen19} or make simplifying assumptions. For example, \cite{GHL21} assumes that the eigenspaces are isolated  and consider the case where lowest eigenvalues are sufficiently well separated from the rest of the spectrum for all parameter values.
\\
\indent  In this work, we design a Ritz method adapted to the particular form of matrix operators involved in problem \eqref{eq:evp}. Contrary to previous methods, we provide error analysis regardless of multiple or crossing eigenvalues.  The subspace is designed by splitting any eigenvector of \eqref{eq:evp} as a sum of component from subspace spanned by eigenvectors of $A_0$ related to the interval of interest and a rest. The rest is related to a correction operator that is an analytic function of $d+1$ arguments. Inspired by \cite{GHL21} and \cite{AndSchw11}, we use sparse stochastic collocation to create a subspace that can approximate range of the correction operator with desired accuracy.
The Ritz space $V$ is then built as a direct sum of two subspaces: one being the span of eigenvectors of $A_0$ related to the interval of interest, and  the second one related to sparse stochastic collocation. The correction operator is similar to the the mapping appearing in the Feshbach-Schur perturbation theory, see \cite{FBS,FBS2} and reference therein.

The structure of the article is as follows.
In the second section, we quickly review the Ritz method and provide an application example of GEVP \eqref{eq:evp}. In section \ref{Section_Gen_eig_beh}, we introduce the notion of average matrix and study  spectral gaps of \eqref{eq:evp}. In section \ref{Section_correction_operator}, we define the correction operator and study its analyticity. Subsection \ref{subsection_Stoch} is dedicated to stochastic collocation and convergence rates of such approximation. In section \ref{Section_approx_av_eig}, we study the accuracy of our Ritz subspace method and provide in a convergence rate in Theorem \ref{Thm_Bauer_Fikes_like}. Finally, in section \ref{num_exp_section}, we illustrate our findings with numerical experiments that validate our theory.


\section{Background}


In this section, we revisit the Ritz method, its error analysis, and give an example of parametric PDE leading to \eqref{eq:evp}. We begin by defining notation used in the following. Assume that eigenvalues of \eqref{eq:evp} are indexed in non-decreasing order and repeated by their multiplicities. This is 
\begin{equation*}
    \lambda_1(\sigma) \leq \lambda_2(\sigma) \leq  \ldots \leq \lambda_n(\sigma) \quad \mbox{for any $\sigma \in S$}. 
\end{equation*}
We denote the corresponding $M$-orthonormal eigenvectors by $( x_k(\sigma))$. Eigenspace related to eigenvalue $\lambda(\sigma)$ and parameter $\sigma \in S$ is denoted by $E_{\lambda(\sigma)}(\sigma)$. The notation $E_{\leq \gamma}(\sigma)$ is used for the direct sum of eigenspaces up to $\gamma$, i.e. 
\begin{equation*}
E_{\leq \gamma}(\sigma) =\oplus_{\lambda_k(\sigma) \leq \gamma} E_{\lambda_k(\sigma)}(\sigma) 
\end{equation*} 
%

Recall that for all $\sigma\in S$, the eigenvalues $\lambda_j(\sigma)$ of \eqref{eq:evp} are given by the famous min-max characterization \cite[Section 10.2]{Parl80}, see also \cite[Chap. 10]{Parl80} or \cite{Johnson85}:
\begin{equation}\label{def_minmax_eig}
    \lambda_j(\sigma)=\min_{U\in V^{(j)}}\max_{x\in U} \frac{x^T A(\sigma)x}{x^T M x}
\end{equation}
for all $j$, where $V^{(j)}=\{ V\subset \mathbb{R}^n, dim(V)= j \}$.
In particular, the extremal eigenvalues satisfy 
\begin{equation}
\label{eq:RQ}
\lambda_1(\sigma)=\min_{x\in \mathbb{R}^n} \frac{x^T A(\sigma) x}
{x^T M x}
\quad \mbox{and}
\quad 
\lambda_n(\sigma)=\max_{x\in \mathbb{R}^n} \frac{x^T A(\sigma) x}
{x^T M x}.
\end{equation}
Further, it is straightforward to show that $\lambda_{k}(\sigma)=\max_{x\in E_{\leq \lambda_k(\sigma)}(\sigma)}\frac{\|x\|_{A(\sigma)}^2}{\|x\|_{M}^2}$ and thus
\begin{equation}\label{Bound_min_max}
    \| x \|_{A(\sigma)}^2 \leq \lambda_k(\sigma) \| x \|^2_{M} \quad\text{for any}\quad x\in E_{\leq \lambda_k(\sigma)}(\sigma).
\end{equation}
Here, and in the following, we denote the inner product and the induced norm defined by $B\in\mathbb{S}^{n\times n}_{++}$ as $\langle .|.\rangle_B$ and $\|.\|_B$ in $\mathbb{R}^n $, respectively. 


\subsection{Review on Ritz method}
\label{subsection_Ritz}


We proceed to describe the Ritz method and derive an upper bound for Ritz eigenvalue error. The material is well known for fixed $\sigma$ and can be found, e.g. in \cite{Saad11}. 

A Ritz method first seeks for a subspace $V$ and then computes an approximate eigensolution of \eqref{eq:evp} from $V$ as follows. Let $\{ v_1,v_2,\ldots,v_k \}$ be a basis of $V$ and denote by $Q\in\mathbb{R}^{n\times k}$ the matrix with column vectors $v_1,v_2,\ldots,v_k$. Then the \emph{Ritz approximation to \eqref{eq:evp} from $V$} is obtained as $(\mu(\sigma),Qy(\sigma))$, where $(\mu(\sigma), y(\sigma)) \in (0,\Lambda) \times \mathbb{R}^k \setminus \{0\}$ satisfy
\begin{equation}
\label{eq:Revp}
    Q^T A(\sigma) Q y(\sigma) = \mu(\sigma) Q^T M Q y(\sigma).
\end{equation}
Assume that Ritz eigenvalues of \eqref{eq:Revp} are indexed in non-decreasing order and repeated by their multiplicities. This is, 
\begin{equation*}
    \mu_1(\sigma) \leq \mu_2(\sigma) \leq \ldots \leq \mu_k(\sigma), \quad \mbox{for $k=\mathop{dim} V$ and any $\sigma \in S$}.
\end{equation*}
The goal of this article is to design a space $V$ so that the Ritz eigenvalue approximation on $(0;\Lambda)$ has the desired accuracy and $\mathop{dim}V$ is small so that GEVP \eqref{eq:Revp}  is inexpensive to solve.

We proceed to study the Ritz eigenvalue error. There exists different upper bounds for the error that all are product of two terms, see \cite{Kynazev,Boffi}. The first term may depend on the applied subspace $V$ as well as the spectrum of the pencil $(A(\sigma),M)$ and can be complicated to estimate. The second term measures how accurately $V$ can approximate either the eigenspace related to the eigenvalue under investigation or eigenspaces related to all previous eigenvalues. In this work, we use a variant where the first term in the bound has a simple dependency on $V$ and the second term measures approximablility of all eigenvectors related to spectral interval of interest $(0;\Lambda)$.  

Before studying Ritz eigenvalue approximation error, we introduce notation and recall technical results. Let $B\in\mathbb{S}_{++}^{n\times n}$.  In the following $P_B$ denotes $B$-orthogonal projection. When orthogonal projection is used, we specify $R(P_B)$ that uniquely defines $P_B$. By $B$-orthogonality, it holds that 
\begin{equation}
\label{eq:P_B_prop1}
    \| x \|^2_B = \| P_B x \|^2_B + \| (I-P_B) x \|_B^2
\end{equation}
for any $x \in \mathbb{R}^n$. Also, note that 
\begin{equation}
\label{eq:best_approximation}
    \| x - P_B x \|_B = \min_{ v \in R(P_B)} \| x - v\|_B.
\end{equation}
This is,  $P_B x$ is the best approximation of $x$ from $R(P_B)$ in $B$-norm.

Next Lemma relates error in the Ritz approximation of $\lambda_k(\sigma)$ to approximability of $E_{\leq \lambda_k(\sigma)}(\sigma)$ in $V$. This result is well known and slightly different versions of it are given in  \cite[Section 4.3]{Saad11} and \cite{Boffi}, see also \cite{Kynazev}.

\begin{Lemma}\label{Lemma_gen_res_bound_lamb_by_x}
Let $\sigma \in S$, $(\lambda_k(\sigma))_k$ be the eigenvalues of \eqref{eq:evp}, and $(\mu_m(\sigma))_m$ their Ritz approximation from $V$. In addition, let $E_{\leq \lambda_k(\sigma)}(\sigma)$ be the union of all eigenspaces related to $(0,\lambda_k(\sigma))$ and $P_{A(\sigma)}$ the $A(\sigma)$ orthogonal projection to $V$. Assume that $e(\lambda_k(\sigma),V) = \max_{x\in E_{\leq \lambda_k(\sigma)}(\sigma)} \|(I-P_{A(\sigma)}) x \|_{M} \|x\|_M^{-1} < 1$. Then there holds that 
\begin{equation}
    \frac{\mu_k(\sigma)-\lambda_k(\sigma)}{\lambda_k(\sigma)} \leq C(\lambda_k(\sigma),V) \max_{x\in E_{\leq \lambda_k(\sigma)}(\sigma), \|x\|_{A(\sigma)}=1} \|(I-P_{A(\sigma)})x \|_{A(\sigma)}^2.
\end{equation}
for any $\sigma \in S$. The multiplicative factor $C(\lambda_k(\sigma), V) = (1-e(\lambda_k(\sigma),V)^2)^{-1}$ approaches one when $e(\lambda_k(\sigma),V)$ tends to zero. 
\end{Lemma}

Before proving this Lemma, we first provide intermediate properties. Indeed, we can see the quantity $\max_{x\in E_{\leq \lambda_k(\sigma)}(\sigma), \|x\|_{A(\sigma)}=1} \|(I-P_{A(\sigma)})x \|_{A(\sigma)}^2 $ as the distance between spaces $E_{\leq \lambda_k(\sigma)}(\sigma)$ and $P_{A(\sigma)}E_{\leq \lambda_k(\sigma)}(\sigma)$.  We start by giving property of distance between two finite dimensional subspaces.

\begin{Lemma}
    Let $H$ be a Hilbert space, $X,Y \subset H$ be finite dimensional subspaces and $\mathop{dim}(X) = \mathop{dim}(Y)$. Then 
    \begin{equation*}
        \dist(X,Y)=\dist(Y,X).
    \end{equation*}
\end{Lemma}
\begin{proof}
    Recall that $dist(X,Y)^2=\sup_{ \substack{x \in X \\ \|x \|_H = 1} } \inf_{ w \in Y } \| x - w\|^2_{H}$. Let $\{x_1,\ldots,x_K\}$ and $\{y_1,\ldots,y_K\}$ be orthonormal basis of $X,Y$ and $P_X$, $P_Y$ $H$-orthogonal projection operators onto $X$ and $Y$ respectively. Then restrictions ${P_X}_{|Y}$ and ${P_Y}_{|X}$ are linear operators with matrix representation $M_{Y\rightarrow X}$ and $M_{X\rightarrow Y}=M_{Y\rightarrow X}^T$, respectively. Then, from orthogonality,
    \begin{equation*}
        \dist(X,Y)^2=
\sup_{ x \in X } \frac{\| (I - {P_Y}_{|X}) x \|^2_{H}}{\| x \|_H^2} = 1 - \min_{x\in X} \frac{\| {P_Y}_{|X} x \|^2_{H}}{\| x \|_H^2}=1-\sigma_{m}(M_{X\rightarrow Y})^2,
    \end{equation*}
    where $\sigma_m(M_{X\rightarrow Y}) \geq 0$ is the smallest singular value of $M_{X\rightarrow Y}$. 
    Similarly, $\dist(Y,X)^2  = 1 - \sigma_{m}(M_{X\rightarrow Y}^T)^2$. Observing that $\sigma_{m}(M^T_{X\rightarrow Y}) = \sigma_{m}(M_{X\rightarrow Y})$ completes the proof.
\end{proof}
\begin{remark}
Assume that $\dim(P_{A(\sigma)}E_{\leq\lambda_k(\sigma)}(\sigma))=\dim(E_{\leq\lambda_k(\sigma)}(\sigma))$. Then it follows from the above Lemma that 
\begin{equation}\label{eq:eq_dist}
       \dist(P_{A(\sigma)}E_{\leq \lambda_{k}(\sigma)}(\sigma),E_{\leq \lambda_{k}(\sigma)}(\sigma))=\dist(E_{\leq \lambda_{k}(\sigma)}(\sigma),P_{A(\sigma)}E_{\leq \lambda_{k}(\sigma)}(\sigma)).
    \end{equation}
In terms of projection matrices,  equality \eqref{eq:eq_dist} is equivalent to
\begin{equation}\label{eq:exchange_PA_PE}
    \|(I-P_E)P_A\|_{A(\sigma)}=\|(I-P_A)P_E\|_{A(\sigma)}.
\end{equation}
where $P_E$ is the $A(\sigma)/M$-orthogonal projection operator onto $E_{\leq \lambda_k(\sigma)}(\sigma)$.
\end{remark}
\begin{proof}[Proof of Lemma \ref{Lemma_gen_res_bound_lamb_by_x}]
    First, we establish that
\begin{equation}
\label{eq:dim}
\mathop{dim} P_{A(\sigma)} E_{\leq \lambda_k(\sigma)}(\sigma) = \mathop{dim} E_{\leq \lambda_k(\sigma)}(\sigma).
\end{equation} 
If $\mathop{dim}P_{A(\sigma)} E_{\leq \lambda_k(\sigma)}(\sigma) < \mathop{dim} E_{\leq \lambda_k(\sigma)}(\sigma)$, there exists $x_1,x_2 \in E_{\leq \lambda_k(\sigma)}(\sigma), x_1 \neq x_2$ s.t. $P_{A(\sigma)}x_1 = P_{A(\sigma)}x_2$. Then $P_{A(\sigma)}(x_1-x_2) = 0$ and $(I-P_{A(\sigma)})(x_1-x_2) = (x_1-x_2)$. Thus $\|(I-P_{A(\sigma)})(x_1-x_2)\|_M \| (x_1-x_2) \|^{-1}_M = 1$ that is a contradiction.
By \eqref{def_minmax_eig} and \eqref{eq:dim}, it holds that
\begin{equation*}
           \mu_k(\sigma)\leq \max_{x\in P_{A(\sigma)} E_{\leq \lambda_k(\sigma)}(\sigma)} \frac{x^T A(\sigma) x}{x^T M x} = \max_{x\in P_{A(\sigma)} E_{\leq \lambda_k(\sigma)}(\sigma)} \frac{\| x\|_{A(\sigma)}^2}{\|  x\|_{M}^2}=\frac{\|\bar{x}\|_{A(\sigma)}^2}{\|\bar{x}\|_M^2}
\end{equation*}
for some $\bar{x}\in P_{A(\sigma)}E_{\leq \lambda_k(\sigma)}(\sigma)$. Define $P_E$ as the $A(\sigma)/M$-orthogonal projection operator onto $E_{\leq \lambda_k(\sigma)}(\sigma)$ and consider the $A(\sigma)/M$-orthogonal decomposition 
\begin{equation*}
    \bar{x}=P_E \bar{x}+(I-P_E)\bar{x}\quad \in \quad  E_{\leq \lambda_k(\sigma)}(\sigma) \oplus E_{\leq \lambda_k(\sigma)}(\sigma)^{\perp}.
\end{equation*}
Then 
\begin{equation}\label{eq:upp_bnd_prf_1}
    \frac{\|\bar{x}\|_{A(\sigma)}^2}{\|\bar{x}\|_M^2}
    \leq \frac{\|P_E\bar{x}\|_{A(\sigma)}^2}{\|P_E\bar{x}\|_M^2}+\frac{\|(I-P_E)\bar{x}\|_{A(\sigma)}^2}{\|\bar{x}\|_M^2} 
    \leq \lambda_k(\sigma)+\frac{\|(I-P_E)\bar{x}\|_{A(\sigma)}^2}{\|\bar{x}\|_M^2}.
\end{equation}
As $\bar{x}\in P_{A(\sigma)}E_{\leq \lambda_k(\sigma)}(\sigma)$, there exists $y\in E_{\leq \lambda_k(\sigma)}(\sigma)$ such that $\bar{x}=P_{A(\sigma)}y$. Moreover, 
\begin{align*}
    \|P_{A(\sigma)}y\|_M^2\geq \Big( \|y\|_M-\|(I-P_{A(\sigma)})y\|_M \Big)^2&=(1-e(\lambda_k(\sigma),V))^2\|y\|_M^2
    \\
    &\geq \frac{(1-e(\lambda_k(\sigma),V))^2}{\lambda_k(\sigma)}\|y\|_{A(\sigma)}^2
\end{align*}
We finally obtain from this and \eqref{eq:upp_bnd_prf_1}
\begin{align*}
    \mu_k(\sigma)&\leq \lambda_k(\sigma)+\frac{\lambda_k(\sigma)}{(1-e(\lambda_k(\sigma),V))^2}\frac{\|(I-P_E)y\|_{A(\sigma)}^2}{\|y\|_{A(\sigma)}^2}\\
     &\leq \lambda_k(\sigma)+\frac{\lambda_k(\sigma)}{(1-e(\lambda_k(\sigma),V))^2}\frac{\|(I-P_E)y\|_{A(\sigma)}^2}{\|P_{A(\sigma)}y\|_{A(\sigma)}^2}\\
     & \leq \lambda_k(\sigma)+\frac{\lambda_k(\sigma)}{(1-e(\lambda_k(\sigma),V))^2} \|(I-P_E)P_{A(\sigma)}\|_{A(\sigma)}^2.
\end{align*}
Using \eqref{eq:exchange_PA_PE}, we obtain the claimed result.
\end{proof}
\begin{Cor}
    \label{RHS}
Make the same assumptions and use the same notation as in Lemma \ref{Lemma_gen_res_bound_lamb_by_x}. In addition, assume that $\lambda_k(\sigma) \in (0,\Lambda)$. Then there holds that 
\begin{equation}\label{eigenvalue_error}
\frac{\mu_k(\sigma)-\lambda_k(\sigma)}{\lambda_k(\sigma)} 
\leq C  \max_{i \in {1,\ldots,N(\sigma) } } \|(I-P_{A(\sigma)})x_i(\sigma) \|_{A(\sigma)}^2, 
\end{equation}
where $C\equiv C(V)\lambda_1(\sigma)^{-1}$ for
\begin{equation*}
C(V) =  \frac{1}{(1-e(\sigma,V))^2}
\quad \mbox{and} \quad 
    e(\sigma,V) = \max_{x\in E_{\leq \Lambda}(\sigma)} \frac{ \|(I-P_{A(\sigma)}) x \|_{M} }{ \|x\|_M}.
\end{equation*}
\end{Cor}

We proceed to outline how the error estimate given in Lemma~\ref{Lemma_gen_res_bound_lamb_by_x} is applied in our context. In the following, we measure the eigenvalue error as 
\begin{equation}\label{eigenvalue_error}
\int_{S} \max_{k\in \{1,\ldots,N(\sigma)\} } \frac{\mu_k(\sigma)-\lambda_k(\sigma)}{\lambda_k(\sigma)} \; d\sigma
\end{equation}
It follows from Corollary~\ref{RHS}  that
\begin{equation}\label{equation_RHS}
\int_S \max_{k\in \{1,\ldots,N(\sigma)\} } \frac{\mu_k(\sigma)-\lambda_k(\sigma)}{\lambda_k(\sigma)} \; d\sigma \leq 
C \int_S \max_{k\in {1,\ldots,N(\sigma)}} \|(I-P_{A(\sigma)})x_k(\sigma) \|_{A(\sigma)}^2 \; d\sigma,
\end{equation}
\noindent where the constant $C$ is as defined in Corollary~\ref{RHS}. Assuming that $e(\sigma,V)<\frac{1}{2}$, $C$ becomes independent on $V$. In the following, we derive an upper bound that considerably simplifies the RHS of \eqref{equation_RHS}.


\subsection{Application in FEM}
\label{subsec_Ex}


\begin{example}
\label{example1}
One motivation for our research comes from solving parameter dependent partial differential equations. 
Let $\Omega = (-1,1)^2$ and $\Omega_1 = B(0,\frac{1}{2})$. Denote the characteristic function of set $\omega \subset \mathbb{R}^d$ as $\chi_\omega$  and define  $\gamma_\sigma = 1 + \chi_{\Omega_1} \sigma$ for any $\sigma\in \mathbb{R}$. Assume that $\sigma \in [-\kappa,\kappa]$ for  some $\kappa < 1$, and consider the parametric eigenvalue problem in it's weak form: find $(\lambda,u) \in \mathbb{R} \times H_0^1(\Omega) \setminus \{0\}$ satisfying 
\begin{equation}
\label{eq:FE_evp}
\int_{\Omega} \gamma_\sigma \nabla u \cdot \nabla v  = \lambda \int_{\Omega} u v \quad \mbox{for all $v \in H^1_0(\Omega)$}. 
\end{equation}
We are interested  in solving \eqref{eq:FE_evp} using the finite element method (see \cite{Braess07,Brenner19} for an introduction in FEM). Let $\{ \phi_i \}_{i=1}^n$ be a basis of a FE-space.  Then the finite element discretization of problem \eqref{eq:FE_evp} is equivalent to the generalized eigenvalue problem: find $(\lambda(\sigma),x(\sigma))\in \mathbb{R} \times \mathbb{R}^n $ satisfying $A(\sigma)x(\sigma)=\lambda(\sigma)Mx(\sigma) $ where
\begin{equation}
    A(\sigma)_{i,j}:=\int_{\Omega}\gamma_\sigma \nabla \phi_i \nabla \phi_j \hspace{0.2cm} \text{ and }\hspace{0.2cm}M_{i,j}:= \int_{\Omega}\phi_i \phi_j \hspace{0.2cm} \text{ for }\hspace{0.2cm}i,j=1,\ldots, n.
\end{equation}
This GEVP is of the type \eqref{eq:evp} with $A(\sigma) = A_0 + A_1 \sigma$. In addition, it holds that $w^T A(\sigma) w = \int_{\Omega} \gamma_\sigma \nabla v_h \cdot \nabla v_h$ and $w^T M w = \int_{\Omega}v^2_h$, where $v_h$ is the FE function corresponding to coefficient vector $w$. It follows from positivity of these bilinear forms that matrices $A(\sigma)$,$M \in \mathbb{S}_{++}^{n\times n}$ for $\gamma_\sigma > 0$ in $ \Omega$.
\end{example}

\begin{example}\label{exemple2.2}
    Problem \eqref{eq:evp}  arises from the stochastic diffusion eigenvalue problem: let $D\subset \mathbb{R}^d$ and $\Omega$ be a suitable probability space as in \cite{Bieri09}. Given  $a\in L^{\infty}(\Omega\times D)$ satisfying $0<a_{\min} \leq a(\omega,x)\leq a_{\max}$ for all $(\omega,x)\in \Omega\times D$, find $(u,\lambda)\in L^2(\Omega\times D)\times \mathbb{R}_{+} $ such that
\begin{equation}\label{model_problem_1}
    \begin{cases}
        -div(a(\omega,\cdot)\nabla u(\omega,\cdot))&=\lambda u(\omega,\cdot),\quad \text{in }D\\
        u(\omega,\cdot)&=0,\quad \text{on }\partial D
    \end{cases}.
    \end{equation}
   The diffusion coefficient function $a$ admits the so-called Karhunen-Loève (KL) expansion
    \begin{equation*}
        a(\omega,x)=\mathbb{E}_{a}(x) + \sum_{m\geq 1} \sqrt{\zeta_m} \varphi_{m}(x)Y_m(\omega) 
    \end{equation*}
    where $\mathbb{E}_{a}(x)=\int_{\Omega}a(\omega,x)dP(\omega)$ and $(\zeta_m,\varphi_m)_{m\geqslant 1}$ is a countable sequence of eigenpairs related to a covariance operator, see \cite{Bieri09}. The $KL$-eigenvalues $\zeta_m$ have decreasing magnitude with $\zeta_1\geq \zeta_2\geq
     \ldots \geq \zeta_m \rightarrow 0 $ as $m\rightarrow \infty$ and the $KL$-eigenfunctions $\varphi_m(x)$ are $L^2$-orthonormal, which all together, and in order to numerically treat the KL-expansion, leads to a truncation $a_M(\omega,x)=\mathbb{E}_{a}(x) + \sum_{m= 1}^M \sqrt{\zeta_m} \varphi_{m}(x)Y_m(\omega)$.
     Using the truncated KL-expansion, considering $Y_M(\omega)$ as  the parameter, and applying FE-method one finally obtains a multi-parametric matrix eigenvalue problem of the form \eqref{eq:evp}.
     \end{example}


\section{Generalized eigenvalue behaviour}
\label{Section_Gen_eig_beh}


In this section, we first recall spectral equivalence between two matrices. Then we review perturbation results that link the spectrum of $A(\sigma) =  A_0 + \sum_{m=1}^d \sigma_m A_m$ to the spectrum of $A_0$. These results follow from the min-max characterization \eqref{def_minmax_eig}. 
Recall that most existing research on multi-parametric eigenvalue problems assume that eigenvalues do not cross in the parameter space. To understand limitations of such approaches, we show how distance between eigenvalues depends on the parameter $\sigma$.


\subsection{Average matrix}


Assume that there exists an average matrix $\overline{A}$ that is spectrally equivalent to $A(\sigma)$ for all $\sigma \in S$. 

\begin{Def}
    Let $\overline{A} \in \mathbb{S}^{n\times n}_{++}$. We say that $ \overline{A} $ is spectrally equivalent to $A:S\rightarrow \mathbb{S}^{n\times n}_{++}$ if there exists $\alpha,\beta \in \mathbb{R}^+$ such that 
\begin{equation}
\label{eq:spequiv}
\alpha v^T \overline{A} v \leq v^T A(\sigma) v \leq \beta v^T \overline{A} v
\end{equation}
for all $v \in \mathbb{R}^n$ and $\sigma \in S$.
\end{Def}
\noindent Note that relevant eigenvalues and eigenvectors of such average matrices can be computed using classical methods. 

\begin{example}\label{example_average_matrix}
The matrix operator $A(\sigma)$ defined in Example~\ref{example1} is spectrally equivalent to the matrix $\overline{A}$ defined by $\overline{A}_{i,j}=\int_{\Omega}\nabla \phi_i \nabla \phi_j $. Indeed, $ \gamma_\sigma \in \begin{bmatrix}
1-\kappa;1+\kappa
\end{bmatrix}$
 for all $\sigma\in [-\kappa;\kappa]$. Then 
 $ v^T A(\sigma)v = \int_{\Omega} \gamma_\sigma \nabla v_h\nabla v_h  \leq (1+\kappa) \int_{\Omega}\nabla v_h \nabla v_h  $ and $ v^T A(\sigma)v = \int_{\Omega} \gamma_\sigma \nabla v_h \nabla v_h  \geq (1-\kappa) \int_{\Omega}\nabla v_h \nabla v_h$. Here $v_h$ is the FE function corresponding to coefficient vector $v$. Hence, $A(\sigma)$ is spectrally equivalent to $\overline{A}$  with $\alpha = 1-\kappa$ and $\beta = 1+\kappa$.
\end{example}

\begin{remark}
Let $A(\sigma) =  A_0 + \sum_{m=1}^d \sigma_m A_m$ and assume that $A_0 \in \mathbb{S}^{n \times n}_{++}$. In this case, the spectral equivalence constants $\alpha,\beta$ are related to the largest possible perturbation  $\sum_{m=1}^d \sigma_m A_m$ of $A_0$ in $S$. 
\end{remark}

 Using \eqref{eq:spequiv}, we can relate eigenvalues of pencil $(A(\sigma),M)$ and those of pencil $(\overline{A},M)$.


\subsection{Spectral separation}


Consider $A(\sigma) =  A_0 + \sum_{m=1}^d \sigma_m A_m$ where $A_0 \in \mathbb{S}^{n \times n}_{++}$. Next we study how large perturbation to $A_0$ is allowed before the eigenvalues may cross. We proceed by linking eigenvalues of $A(\sigma)$ to those of $A_0$ using spectral equivalence as in \cite{GHL21}. These bounds are then used to estimate the gap between eigenvalues of $A(\sigma)$. 

Recall that for all $\sigma\in S$, the eigenvalues $\lambda_j(\sigma)$ of \eqref{eq:evp} are given by the min-max characterization \eqref{def_minmax_eig}.
In particular, 
\begin{equation*}
    \lambda_j(0)=\min_{U\in V^{(j)}}\max_{x\in U} \frac{x^T A_0 x}{x^T M x}.
\end{equation*}
Assume that $A(\sigma)$ in \eqref{eq:evp} is spectrally equivalent to $A_0=A(0)=:\overline{A}\in\mathbb{S}^{n\times n}_{++}$ with parameters $\alpha$ and $\beta$. It follows from the min-max characterization of eigenvalues that
\begin{equation}\label{eq_encad_overA}
		\alpha\lambda_{i}(0)\leq \lambda_{i}(\sigma)\leq \beta\lambda_{i}(0).
		\end{equation}
  for any $\sigma \in S$ and $i \in \{1,\ldots,n\}$.
 That is, equation \eqref{eq_encad_overA} states that eigenvalues of the pencil $(A_0,M)$ and spectral equivalence constants $\alpha,\beta$ give two sided bounds for eigenvalues of \eqref{eq:evp}. We proceed to study how the relative distance between eigenvalues of \eqref{eq:evp} depends on spectral equivalence constants $\alpha,\beta$. Previous works on multiparametric eigenvalue problems assume that this  perturbation is so small that the gap between eigenvalues stays positive for any $\sigma \in S$. 
\begin{Lemma}\label{lemma_isolated_eigspa}
	Let $i \in \{1,\ldots,n\}$ and define $ \delta_0 = \frac{\lambda_{i+1}(0)-\lambda_{i}(0)}{\lambda_{i}(0)} $. If $\alpha$ and $\beta$ are such that $\delta_0>\frac{1}{\alpha}(\beta-\alpha)$, 
then for all $\sigma\in S$, $\frac{\lambda_{i+1}(\sigma)-\lambda_{i}(\sigma)}{\lambda_{i}(\sigma)} \geq \delta $ with
	\begin{equation*}
	\delta=\frac{1}{\beta}\Big(\alpha\delta_0 +(\alpha-\beta)\Big)>0.
	\end{equation*}
\end{Lemma}
A similar result is given in \cite{GHL21} and it tells us that if the $i^{th}$ and $(i+1)^{th}$ eigenvalues of $A_0$ are well separated and $\alpha,\beta$ are sufficiently small, then  the $i^{th}$ and $(i+1)^{th}$ eigenvalues of $A(\sigma)$ do not cross.

\section{Correction operator and the stochastic collocation method}
\label{Section_correction_operator}


The aim of this work is to provide a new way to define a Ritz subspace for multiparametric eigenvalue approximation. Recall, that by \eqref{equation_RHS} Ritz eigenvalue error on $(0,\Lambda)$ is bounded by approximability of corresponding eigenvectors in the Ritz subspace. In this section, we first show how to split the eigenvectors into two components. The first component lies in subspace containing eigenvectors of an average matrix related to eigenvalues in $ [0;\rho\Lambda] $ for some $\rho>1$. The second component is expressed by a correction operator, that is an analytic function of $d+1$ arguments. We then build our Ritz subspace from eigenvectors of the average matrix and samples of the correction operator. The samples are evaluated at interpolation points related to the stochastic collocation operator from \cite{AndSchw11}. This allows us to bound Ritz eigenvalue error on $(0,\Lambda)$ by interpolation error estimates. 

Previously such work has been done in an attempt to compute eigenvalues of the self-adjoint Schrödinger operators. In \cite{Stamm23}, these eigenvalues are then approximated by eigenvalues solutions of a corrected operator on a subspace spanned by the eigenfunctions of the Laplacian operator. However, in our approach, we avoid using a truncation of Neumann series as done in \cite{Stamm23}.


\subsection{Correction formula}
\label{subsec_correction_formula}


Assume that $A : S\rightarrow \mathbb{S}^{n\times n}_{++} $ is spectrally equivalent to an average matrix $\overline{A} \in \mathbb{S}^{n\times n}_{++}$. Let $ \rho\Lambda $ be the length of the sampling interval and recall that
\begin{equation*}
    \overline{E}_{\leq \rho\Lambda}=\{v\in \mathbb{R}^n \hspace{0.2cm}  | \hspace{0.2cm} \exists y\in(0;\rho\Lambda) \text{ s.t. } \overline{A}v=yMv \}
\end{equation*}
this is, $\overline{E}_{\leq\rho\Lambda}$ is the union of all eigenspaces of the pencil $(\overline{A},M)$ associated to eigenvalues on $ (0;\rho\Lambda)$. Let $W$ be a basis matrix of $\overline{E}_{\leq\rho\Lambda}$.  The space $\overline{E}_{\leq\rho \Lambda}$ and its $\overline{A}$-orthocomplement 
\begin{equation*}
    \overline{E}_{\leq\rho\Lambda}^\perp := \{ v \in \mathbb{R}^n \; | \; v^T \overline{A} W = 0 \; \}
\end{equation*}
are central in the rest of this paper. It holds that
\begin{equation}
w_\perp^T \overline{A} w_\perp \geq \rho \Lambda w_\perp^T M w_\perp \quad \mbox{for all $w_\perp \in \overline{E}_{\leq\rho\Lambda}^\perp$}
\end{equation}
Let $W_\perp$ be the basis matrix of $\overline{E}_{\leq\rho\Lambda}^\perp$. Note that $W$ can be obtained by solving few lowest eigenmodes of the pencil $(\overline{A},M)$, whereas $W_{\perp}$ is large, dense, and never available in a practical computation. In the following, we use the notation 
\begin{equation}
B_\perp:=W_\perp^T B W_\perp, \hspace{0.3cm} \forall B\in\mathbb{R}^{n\times n}.
\end{equation}

We proceed to split eigenvectors of \eqref{eq:evp} as the sum of a component in $\overline{E}_{\leq\rho\Lambda}$ and a correction term in $\overline{E}_{\leq\rho\Lambda}^\perp$. A suitable splitting is given in Lemma \ref{lemma_def_Z} after a technical result. 
\begin{Lemma}\label{lemma_anal_Z}
Assume that $ A:S\rightarrow \mathbb{S}^{n\times n}_{++} $ is spectrally equivalent to $\overline{A}$ as in \eqref{eq:spequiv} with constants $\alpha$ and $\beta$. Let $\Lambda > 0$ and $\rho$ be such that $\alpha \rho>1$ and $\rho>1$. Then the matrix $A_{\perp}(\sigma)-tM_{\perp}$ is positive definite for all $(\sigma,t)\in S\times(0; \Lambda)$. 
\end{Lemma}

\begin{proof}
    For all $(\sigma,t)\in S\times(0; \Lambda)$ and $w_\perp\in \overline{E}_{\leq \rho \Lambda}^\perp\backslash \{0\}$, from the spectral equivalence of $A(\sigma)$ with $\overline{A}$, we get
    \begin{equation}
        w_\perp^T A(\sigma)w_\perp -tw_\perp M w_\perp\geq \alpha w_\perp^T \overline{A}w_\perp -tw_\perp M w_\perp\geq (\alpha \rho \Lambda-t)\|w_\perp\|^2_M>0
    \end{equation}
    as $w_\perp\in \overline{E}_{\leq\rho\Lambda}^\perp$. 
\end{proof}

\begin{Lemma} \label{lemma_def_Z} Make the same assumptions as in Lemma \ref{lemma_anal_Z}. Let $\sigma \in S$ and $(\lambda(\sigma),x(\sigma))\in (0;\Lambda) \times \mathbb{R}^n \setminus \{0\}$ be a solution to \eqref{eq:evp}. In addition, let $W_\perp$ be a basis matrix of $\overline{E}_{\leq\rho\Lambda}^\perp$. Then there holds that  
\begin{equation}
\label{eq:corr}
x(\sigma) 
= 
\overline{x}(\sigma) 
+ 
Z(\sigma,\lambda(\sigma)) \overline{x}(\sigma), 
\end{equation}
for $\overline{x}(\sigma) \in \overline{E}_{\leq\rho\Lambda}$ and correction operator $Z(\sigma, t ) : S \times (0,\Lambda) \mapsto \mathbb{R}^{n\times n}$ defined as 
\begin{equation}\label{first_def_Z}
    Z(\sigma,t )=W_\perp(W_\perp^T (A(\sigma)-tM)W_\perp)^{-1}W_\perp^T\delta A(\sigma)
\end{equation}
 with $\delta A(\sigma)=A(\sigma)-\overline{A}$. 
\end{Lemma}

\begin{proof}
    Let $\sigma\in S$ and split $ x(\sigma)=\overline{x}(\sigma)+r(\sigma)$ with $\overline{x}(\sigma) \in \overline{E}_{\leq\rho\Lambda}$ and $r(\sigma)=W_{\perp}\gamma(\sigma)\in \overline{E}_{\leq\rho \Lambda}^\perp$. Then as $(\lambda(\sigma),x(\sigma))\in (0;\Lambda)\times \mathbb{R}^n$ is a solution of \eqref{eq:evp}, we get
    \begin{equation}
        (A(\sigma)-\lambda(\sigma)M)\overline{x}(\sigma)=-(A(\sigma)-\lambda(\sigma)M)W_\perp \gamma(\sigma).
    \end{equation}
    By orthogonality, we first observe that $W_{\perp}^T M\overline{x}(\sigma)=0$ and $W_{\perp}^T \overline{A}\overline{x}(\sigma)=0$.
    By multiplication with basis matrix $W_{\perp}^{T}$ and using these identities gives
    \begin{equation}\label{eq:W_equa}
        W_\perp^T (A(\sigma)-\overline{A})\overline{x}(\sigma) = -W_{\perp}^{T}(A(\sigma)-\lambda(\sigma) M)W_\perp \gamma(\sigma). 
    \end{equation}
Then from Lemma \ref{lemma_anal_Z}, matrix in the right-hand side of equation \eqref{eq:W_equa} is invertible and we obtain
    \begin{align*}
        -W_\perp (W_\perp^T (A(\sigma)-\lambda(\sigma)M)W_\perp)^{-1}W_\perp^T \delta A(\sigma)\overline{x}(\sigma)=W_\perp \gamma(\sigma).
    \end{align*}
\end{proof}

It is straightforward to estimate the $\overline{A}$-norm of  $\overline{x}(\sigma)$ appearing in \eqref{eq:corr}. Recall that $\overline{x}(\sigma)^T \overline{A} W_\perp = 0$ and using spectral equivalence \eqref{eq:spequiv} gives 
\begin{equation}
\| \overline{x}(\sigma) \|_{\overline{A}} \leq \| x(\sigma)\|_{\overline{A}} \leq \frac{1}{\alpha} \| x(\sigma) \|_{A(\sigma)}.
\end{equation}
As $x(\sigma)$ is an eigenvector of the pencil $(A(\sigma),M)$ corresponding to an eigenvalue on $(0,\Lambda)$ and $\| x(\sigma) \|_M = 1$, it follows that 
\begin{equation}
\| x(\sigma) \|_{A(\sigma)} \leq \Lambda^{1/2} \quad \mbox{and thus} \quad \| \overline{x}(\sigma) \|_{\overline{A}} \leq \frac{1}{\alpha} \Lambda^{1/2}. 
\end{equation} 
These estimates are used later.

In the following, we need to evaluate the action of $Z$. Recall that this requires the solution of a linear system that is posed in the subspace $\overline{E}^\perp_{\leq\rho \Lambda}(\overline{A},M)$. This is a non-trivial task, as the matrix $\Wperp$ is not available in any practical computation. We solve this issue in the following Lemma by applying a saddle-point formulation.

\begin{Lemma}\label{system_value_Z}
Let $(\sigma,t)\in S\times (0;\Lambda)$. Make the same assumptions as in Lemma \ref{lemma_anal_Z}. Let $W,W_{\perp}$ be the basis matrices of $\overline{E}_{\leq\rho\Lambda}$ and $\overline{E}_{\leq\rho \Lambda}^\perp$, respectively. In addition, let $K(\sigma,t):=A(\sigma) -t M$, $b \in \mathbb{R}^n$, and $r$ satisfy 
\begin{equation}\label{mat_syst_solv}
\begin{bmatrix} K(\sigma,t) & \overline{A} W \\ W^T \overline{A} & 0 \end{bmatrix} \begin{bmatrix} r \\ \eta \end{bmatrix} = \begin{bmatrix}
\delta A(\sigma) b \\ 0 \end{bmatrix}.
\end{equation}
Then $r = Z(\sigma,t) b$.
\end{Lemma}
\begin{proof}
Second line of \eqref{mat_syst_solv} gives that $W^{T}\overline{A}r=0$, i.e. $r \in \overline{E}_{\leq\rho\Lambda}^{\perp}$. 
Then the first line \eqref{mat_syst_solv} with $r=W_\perp\alpha$ yields
\begin{equation}
    (A(\sigma)-tM)W_\perp \alpha +\overline{A}W\eta=\delta A(\sigma) b.
\end{equation}
Considering that $W_{\perp}^{T}\overline{A}W=0$ by orthogonality, we obtain
\begin{equation}
    W_\perp^T(A(\sigma)-tM)W_\perp \alpha =W_\perp^T \delta A(\sigma) b
\end{equation}
and then
\begin{equation}
    W_\perp \alpha = W_\perp(W_\perp^T (A(\sigma)-tM)W_\perp)^{-1}W_\perp^T \delta A(\sigma) b.
\end{equation}
This completes the proof.
\end{proof}

We proceed to prove that $Z$ admits complex analytic extension to the complex plane with respect to the parameter vector $\sigma$. This result will allow to approximate $Z$ on $\sigma$ by interpolation. First we establish that there exists a neighbourhood of $S\times(0,\Lambda)$ in which $Z$ is well defined. 
\begin{Lemma}\label{lemma:Cinv} Assume that $ A:S\rightarrow \mathbb{S}^{n\times n}_{++} $ is spectrally equivalent to $\overline{A}$ as in \eqref{eq:spequiv} with constants $\alpha$ and $\beta$. Decompose $\sigma \in \mathbb{C}^d$ as $\sigma = \sigma_S + \sigma_P$, where $\sigma_S \in S$ is the closest point to $\sigma$ from $S$. Let $\Lambda > 0$ and $\rho$ such that $\alpha \rho>1$ and $\rho>1$. Then 
$N(A_\perp(\sigma)-tM_\perp) = \{0\}$ for any $t \in (0,\Lambda)$ and $\sigma \in D$, where 
\begin{equation*}
    D := \{ \sigma \in \mathbb{C}^d \; | \; \max_{x\in \mathbb{C}^{m}} \frac{| x^*\delta A_\perp(\sigma_P)x|}{x^* M_{\perp} x} < (\alpha\rho-1)\Lambda \}.
\end{equation*}
Here $m = \mathop{dim} E_{\leq \rho \Lambda}^{\perp}$ and $\delta A_{\perp}(\sigma_P):=A_{\perp}(\sigma_P)-A_{0\perp}$. 
\end{Lemma} 
\begin{proof}
Let $z\in \mathbb{C}^n$ be such that $\left( A_\perp(\sigma)-t M_\perp \right) z = 0$. The matrix $A_\perp(\sigma)-tM_\perp$ has a trivial null-space  only if $z = 0$. By linearity 
\begin{equation*}
\left( A_\perp(\sigma)-t M_\perp \right) z 
 = \left( A_\perp(\sigma_S)-t M_\perp \right) z + \delta A_\perp(\sigma_P) z 
\end{equation*}
Multiplying with $z^* \in \mathbb{C}^n$, taking absolute value and using the inverse triangle inequality gives
\begin{equation*}
0 = |z^* \left( A_\perp(\sigma)-t M_\perp \right) z | \geq 
\left|\; | z^* \left( A_\perp(\sigma_S)-t M_\perp \right) z | - | z^* \delta A_\perp(\sigma_P) z | \; \right|. 
\end{equation*}
We proceed to estimate the first term on the RHS from below. As $(A_\perp(\sigma_S) - tM_\perp)$ is real s.p.d. it holds that
\begin{equation*}
\frac{ z^* \left( A_\perp(\sigma_S) - tM_\perp \right)z }{z^* M_{\perp} z } \geq \min_{y \in \mathbb{R}^m} \frac{ y^T A_\perp(\sigma_S) y - t y^T M_\perp y}{y^T M_\perp y}.
\end{equation*}
By spectral equivalence, it holds that 
\begin{equation*}
\frac{ z^* \left( A_\perp(\sigma_S) - tM_\perp \right)z }{z^* M_{\perp} z } \geq (\alpha \rho - 1) \Lambda
\end{equation*}

If $\sigma \in D$, we then obtain $0 \geq C \| z \|_M^2$ for some $C>0$. Hence, $z = 0$ that completes the proof. 

\end{proof}
\begin{Cor}\label{cor:1}
Make the same assumptions and use the same notation as in Lemma~\ref{lemma:Cinv}. Then $Z:D \mapsto \mathbb{C}^n$ is complex-analytic for any $t\in(0,\Lambda)$.
\end{Cor}
\begin{proof} We proceed by showing that $(A_\perp(\sigma) - tM)^{-1}$ admits power series representation with respect to $\sigma-\sigma_0$ in the neighborhood of $\sigma_0\in D$ and $t\in(0,\Lambda)$. By Lemma \ref{lemma:Cinv}, $(A_\perp(\sigma_0) - tM_\perp)$ has an inverse. Hence, 
\begin{equation*}
    (A_\perp(\sigma) - tM_\perp) = (A_\perp(\sigma_0) - tM_\perp) \left( I+(A_\perp(\sigma_0) - tM_\perp)^{-1} \delta A_{\perp}(\sigma-\sigma_0) \right).
\end{equation*}
For sufficiently small $\| \sigma - \sigma_0 \|$, it holds that $\| (A_\perp(\sigma_0) - tM_\perp)^{-1} \delta A_{\perp}(\sigma-\sigma_0) \| < 1$. Hence, the above matrix is invertible and has the series expansion
\begin{equation*}
(A_\perp(\sigma) - tM_\perp)^{-1} = \sum_{i=0}^{\infty} \left( K_\perp(\sigma_0,t)^{-1} \delta A_{\perp}(\sigma-\sigma_0) \right)^i K_\perp(\sigma_0,t)^{-1},
\end{equation*}
where $K_\perp(\sigma_0,t) = (A_\perp(\sigma_0) - tM_\perp)$. Thus, $Z(\sigma,t)$ has a multiparametric power series expansion with respect to $\sigma-\sigma_0$ in the neigbourhood of any $\sigma_0 \in D$. This completes the proof. 
\end{proof}

\begin{example}\label{example_2} 
     Solution of multiparametric eigenvalue problems has previously been studied in \cite{GHL21} under the assumption that $A$ is defined by $A(\sigma)=A_0+\sum_{m=1}^{d}\sigma_m A_m $ with $A_0\in \mathbb{S}^{n\times n}_{++}$ and $A_m\in \mathbb{S}^{n\times n}$ satisfying  $\|A_0^{-1/2}A_m A_0^{-1/2}\|_2\leq \kappa_m$ for $m=0,\ldots,d$ and $\| \kappa \|_1 < 1$. These assumptions are stronger than spectral equivalence, as it immediately follows that
    \begin{equation}
        (1-\|\kappa\|_{1})x^{T}A_0 x\leq x^{T}A(\sigma)x\leq (1+\|\kappa\|_{1})x^{T}A_0 x,\quad \forall \sigma\in [-1;1]^{\infty}, \forall x\in \mathbb{R}^n.
    \end{equation}
    Under these conditions, the 
    correction operator $Z$ can be extended analytically to an operator on $E(\tau)=\{ z\in\mathbb{C}^{d}\quad | \quad |z_m|< \tau_m=\tau_m(\kappa) \}$ where $\tau_0=\frac{(1+\varepsilon)\rho\Lambda}{2}$, $\tau_m=\frac{1-\epsilon}{2 \kappa_m\|\kappa\|_{1}}$ for all $m= 1,\ldots,d$ and $\varepsilon\in ]0;1[$. 
    In fact if $z=(t,\zeta_1,\zeta_2,\ldots,\zeta_d)\in E(\tau)$,
    define  $X(\zeta,t)=A_{0\perp}^{-1/2}(\sum_{m=1}^{d}\zeta_m A_{m\perp} -tM_{\perp} )A_{0\perp}^{-1/2}$. Then 
    \begin{align*}
        \|X\|_2^2&\leq \sum_{m=1}^{d}|\zeta_m| \|A_{0\perp}^{-1/2}A_{m\perp}A_{0\perp}^{-1/2}\|_2^2 +|t|\|A_{0\perp}^{-1/2} M_{\perp}A_{0\perp}^{-1/2}\|_2^2 \\
        & \leq \sum_{m=1}^{d}|\zeta_m| \kappa_m^2 +\frac{|t|}{\rho\Lambda}< \frac{(1-\varepsilon)}{2}\sum_{m=1}^{d} \frac{\kappa_m}{\|\kappa\|_{1}} +\frac{(1+\varepsilon)\rho\Lambda}{2\rho\Lambda} =\frac{(1-\varepsilon)}{2} +\frac{(1+\varepsilon)}{2}=1 
    \end{align*}
    and so from Neumann's Lemma, the operator $Z(\zeta,t)=W_{\perp}A_{0\perp}^{-1/2}(I+X(\zeta,t))^{-1}A_{0\perp}^{-1/2}W_{\perp}^{T}\delta A(\zeta)$ is well defined in $E(\tau)$ and by identical arguments as in Corollary~\ref{cor:1} also complex analytic in $E(\tau)$.
    \end{example}

We conclude the section by giving two technical Lemmas that characterize $Zy$ for $y\in \overline{E}_{\leq\rho \Lambda}$. These results are useful later when estimating Ritz eigenvalue error. Both Lemmas are based on representing the term $(A_\perp(\sigma) - t M_\perp)^{-1}$, appearing in the definition of $Z$, by utilizing eigenbasis of the pencil $(A_\perp(\sigma),M_{\perp})$. 

 \begin{Lemma}\label{Lemma_express_Z}
     Let $\sigma\in S\subset \mathbb{R}^d$,  $t\in(0;\Lambda)$, and $A:S\rightarrow \mathbb{S}^{n\times n}_{++}$. Assume that $A$ is spectrally equivalent to $\overline{A}\in\mathbb{S}^{n\times n}_{++}$ with constants $\alpha$ and $\beta$ as in \eqref{eq:spequiv}. Then for all $y\in \overline{E}_{\leq \rho\Lambda}$, it holds that
     \begin{equation}
         Z(\sigma,t)y=W_{\perp}\sum_k u_k u_k^T\frac{1}{\xi_k(\sigma)-t}W_{\perp}^TA(\sigma)y
     \end{equation}
     where $\xi_k(\sigma) \geq \rho\Lambda$ and $u_k\in \overline{E}_{\leq\rho\Lambda}^{\perp}$ are such that for all $k$, 
    $ A_{\perp}(\sigma)u_k=\xi_k(\sigma) M_{\perp}u_k$ and $(u_k)_k$ is $M$-orthonormal. 
 \end{Lemma}
 \begin{proof}
     Let $(\xi_k(\sigma),u_k) \in (\rho\Lambda,\infty) \times \mathbb{R}^{\mathop{dim}(\overline{E}^\perp_{\leq \rho \Lambda})}$ be such that for all $k$, 
    \begin{equation*}
        A_{\perp}(\sigma)u_k=\xi_k(\sigma) M_{\perp}u_k \quad \| u_k \|_M = 1.
    \end{equation*}
    Then it holds that $M_{\perp}^{-1/2} A_\perp(\sigma) M_{\perp}^{-1/2} = \sum_{k} \xi_k(\sigma) \hat{u}_k \hat{u}^T_k$ for orthonormal $\hat{u}_k = M_\perp^{1/2} u_k$. Hence,    
    \begin{equation*}
(A_\perp(\sigma) - t M_{\perp})^{-1} = M_\perp^{-1/2} ( M_{\perp}^{-1/2} A_\perp(\sigma) M_{\perp}^{-1/2} - t)^{-1} M_\perp^{-1/2} = \sum u_k u_k^T \frac{1}{\xi_k(\sigma) -t}  
\end{equation*}
    and we obtain the claimed result using definition of $Z$ and observing that $W_\perp^T \overline{A} y = 0$. 
 \end{proof}

 \begin{Lemma}\label{Lemma_bound_Zy}
     Let $\sigma\in S\subset \mathbb{R}^d$,  $t\in(0;\Lambda)$, and $A:S\rightarrow \mathbb{S}^{n\times n}_{++}$. Assume that $A$ is spectrally equivalent to $\overline{A}\in\mathbb{S}^{n\times n}_{++}$ with constants $\alpha$ and $\beta$ as in \eqref{eq:spequiv}. In addition, let $\rho>1$ be such that $\rho\alpha>1$. It holds that for all $y\in \overline{E}_{\leq \rho\Lambda}$
     \begin{equation}
         \|Z(\sigma,t)y\|_{A(\sigma)}
         \leq  \frac{\rho}{(\rho-1)}(\beta\rho\Lambda)^{1/2} \|y\|_{M}.
     \end{equation}
 \end{Lemma}
 \begin{proof}
     Let $(\xi_k(\sigma))_k$ and $(u_k)_k$ be as in Lemma \ref{Lemma_express_Z}. In addition, let $r=\sum_k \gamma_k u_k \in\mathbb{R}^{\mathop{dim}(\overline{E}^\perp_{\leq \rho \Lambda})}$ be such that $(A(\sigma)_{\perp}-tM_{\perp})r=W_{\perp}^T A(\sigma )y,
$ so that $Zy=W_{\perp}r$. Then 
\begin{equation*}
    r^T (A(\sigma)_\perp-t M_\perp) r = \sum_{k} \left(1 - \frac{t}{\xi_k(\sigma)} \right) \xi_k(\sigma) \gamma_{k}^2 \geq \left(1-\frac{1}{\rho} \right) \sum_k \xi_k(\sigma)\gamma_k^2 = \left(1-\frac{1}{\rho} \right) \| W_{\perp} r \|^2_{A(\sigma)}.
\end{equation*}
Thus, $\left(1-\frac{1}{\rho} \right)\| W_{\perp} r \|^2_{A(\sigma)}\leq(W_{\perp}r)^TA(\sigma)y\leq \|W_{\perp}r\|_{A(\sigma)}\|y\|_{A(\sigma)}$, and hence, $\left(1-\frac{1}{\rho} \right) \| W_{\perp} r \|_{A(\sigma)}\leq \| y\|_{A(\sigma)}$. Further, by spectral equivalence, definition of $\overline{E}_{\leq \rho \Lambda}$ and \eqref{Bound_min_max}, $\| y \|_{A(\sigma)}\leq (\beta\rho\Lambda)^{1/2} \|  y\|_{M}$ and then
\begin{equation}
\|Z(\sigma,t)y\|_{A(\sigma)} = \|W_{\perp}r\|_{A(\sigma)}\leq \frac{\rho}{(\rho-1)}(\beta\rho\Lambda)^{1/2} \|y\|_{M}.
\end{equation}
 \end{proof}


\subsection{Ritz subspace by polynomial approximation}
\label{subsection_Stoch}


Now that we have all the properties from the previous section, we return to building our Ritz subspace $V$. In this work, we define $V$ as
\begin{equation}
    V=\overline{E}_{\leq\rho\Lambda}\oplus X_{\rho\Lambda}
\end{equation}
where the subspace $X_{\rho\Lambda}\subset \overline{E}_{\leq\rho\Lambda}^{\perp}$ is chosen to contain approximations of $Z(\sigma,\lambda(\sigma))\bar{x}(\sigma)$. 

In this section we build $X_{\rho \Lambda}$ related to a polynomial interpolation of the correction operator using collocation points $\{\sigma_i\}\subseteq S$, $\{t_j\}\subseteq (0;\Lambda)$. Two families of collocation points are used as $L^\infty(0,\Lambda)$-convergence with respect to $t$ variable is needed later whereas $L^2(S)$-convergence is sufficient for $\sigma$.  Recall that $W=[w_1,\ldots,w_m]$ is the basis matrix of $\overline{E}_{\leq\rho\Lambda}$. Define sample vectors $s_{i,j,k}$ as 
\begin{equation*}
s_{i,j,k}= Z(\sigma_i,t_j)w_k
\end{equation*} 
and the subspace $X_{\rho\Lambda}$ as 
\begin{equation*}
X_{\rho\Lambda} :=span_{i,j,k} s_{i,j,k}.
\end{equation*} 
The rational for this choice is that such $X_{\rho \Lambda}$ contains the polynomial interpolation of the correction operator, i.e., 
\begin{equation}\label{def_Lag}
    \sum_{i,j} Z(\sigma_i,t_j) \overline{x}(\sigma) \ell_i(\sigma) \ell_j(t) \in X_{\rho \Lambda}\quad \text{for any}\hspace{0.2cm} \sigma\in S\text{ and }t\in (0;\Lambda)
\end{equation}
In this case, we prove an upper bound to the error between solutions of \eqref{eq:evp} and their Ritz approximation in $V$ by estimating how accurately \eqref{def_Lag} approximates $Z(\sigma,t)\overline{x}(\sigma)$ on $S\times (0;\Lambda)$. In the rest of this section, we review the interpolation method and error estimate given in \cite{BNT07} with specific interpolation points which will be used to build $X_{\rho\Lambda}$. 


\subsubsection{The model problem and hierarchic polynomial approximation}

     
    In Example \ref{exemple2.2}, we obtain a solution $u_M$ of the parametric deterministic problem that can be expressed in terms of multivariate orthogonal polynomials and one can consider a truncation of its expansion in this family. However, instead of using a truncation after $M$ terms, one could use the idea of M. Bieri, R. Andreev and C. Schwab \cite{BieriSchwab09} of finding an index set corresponding to the 'active' terms which realizes 'quasi-best N-terms approximation' and will be identified in an index set $\Lambda\subset (\mathbb{N}^{\infty})_c$. 
    
    Once such set $\Lambda$ is given, we need to choose a polynomial basis. Considering that the tensorized Legendre polynomials $\mathcal{L}_{\alpha}=\mathcal{L}_{\alpha_1}\mathcal{L}_{\alpha_2}\ldots$ form an orthogonal family on each variable, spanning $L^2(\Omega)$, we mix this family with the Smoljak method as used in \cite{GHL21}. Indeed in \cite{Smoljak63}, S.A. Smoljak was looking for a new accurate quadrature formula expressed in terms of tensor products in the form $(\sum_{p=1}^{N_1}\gamma_p^{(1)}\tau_p^{(1)})\otimes \ldots \otimes (\sum_{p=1}^{N_s}\gamma_p^{(s)}\tau_p^{(s)}) $ to certain classes of functions. 
    Using tensor products of differences $\nu_v -\nu_{v-1}$ with $\nu_v$ a linear combination of interpolation values of the function, he obtained a quadrature formula that, up to a logarithm in N the number of nodes used, gives the best possible accuracy among all interpolation formulas for the same nodes. 
    Thanks to all these results, we are now able to define the stochastic collocation operator that will be used to approximate solutions with great accuracy.

\subsubsection{The Stochastic collocation operator}
\label{subsubsection_stoch_colloc_def}
 We quickly review some notation from \cite{AndSchw11}: 
 for a sequence $ \eta=(\eta_m)_m $ such that $ \eta_1\geq \eta_2\geq \ldots $ and $ \eta_m\rightarrow 0 $ as $ m\rightarrow \infty $ and a fixed $ \varepsilon>0 $, define 
\begin{equation}\label{def_Aeps}
    \Lambda_{\varepsilon}(\eta)=\{ \alpha\in (\mathbb{N}^\infty)_c \hspace{0.2cm}| \hspace{0.2cm} \eta^\alpha  \geq \varepsilon \}
\end{equation}
where $ \eta^\alpha=\prod_{m\geq 1}\eta_m^{\alpha_m}$ and $ \mathbb{N}^\infty_c=\{ \alpha\in \mathbb{N}^\infty \hspace{0.2cm} | \hspace{0.2cm} supp(\alpha)<\infty \} $ when $ supp(\alpha)=\{ m\geq 1, \alpha_m\neq 0 \} $. 

\begin{Def}[The stochastic collocation operator]
    For $p\in\mathbb{N}$, let $ \{ \chi_{i}^{(p)} \}_{i=0}^{p}\subset [-1;1] $ be the abscissae of the univariate Legendre polynomial of degree $p$ and $ \{ \ell_{i}^{(p)} \}_{i=0}^p $ be the standard Lagrange basis polynomials. For any function $f\in C([-1;1])$, define 
    \begin{equation}\label{def_I_p}
    I_pf(x)=\sum_{j=0}^{p}f(\chi_j^{(p)})\ell_j^{(p)}(x)
    \end{equation}
    the Lagrange polynomial which interpolates $f$ at points $\{ \chi_{j}^{(p)} \}_{j=0}^{p}$. For a decreasing sequence $\eta=(\eta_1,\eta_2,\ldots)$ such that $1>\eta_1\geq \eta_2\geq \ldots$ and $\eta_j\rightarrow 0$ as $j\rightarrow \infty$ and for $\varepsilon>0$, consider $\Lambda_{\varepsilon}(\eta)=\{ \alpha\in (\mathbb{N}^{\infty})_c | \prod_{j\geq 1}\eta_j^{\alpha_j}>\varepsilon \}$ and define {\it the sparse collocation operator}
    \begin{equation}\label{def_inf_Stoch_op}
        \mathcal{I}_{\Lambda_{\varepsilon}(\eta)}:=\sum_{\alpha\in \Lambda_{\varepsilon}(\eta)}\bigotimes_{m\in supp(\alpha)}(I_{\alpha_m}^{(m)}-I_{\alpha_{m}-1}^{(m)})
    \end{equation} 
    where $I_{\alpha_m}^{(m)}$ denotes the operator \eqref{def_I_p} applied on the $m^{th}$ variable in the Kronecker product \eqref{def_inf_Stoch_op}.
\end{Def}
\begin{remark}
As $\eta_m  \rightarrow 0$ for any $\alpha\in \Lambda_{\varepsilon}(\eta)$, there is only a finite number of indices in $supp(\alpha)$ and then $\mathcal{I}_{\Lambda_\varepsilon(\eta)}$ is well defined on bounded continuous functions on bounded domain.
 Note that here for the case when $ A(\sigma)=A_0+\sigma A_1 $, $\eta_1\in \mathbb{R}$ and since $\varepsilon>0$ is fixed, $\Lambda_\varepsilon(\eta)$ is finite.
\end{remark}

\begin{remark}\label{rem_I_Lambda_d}
    Note that a Stochastic collocation operator can be designed for any finite number of variables as $ \mathcal{I}_{\Lambda_{\varepsilon}(\eta)}=\mathcal{I}_{\Lambda_\varepsilon(\eta_1,\ldots,\eta_d)}:=\sum_{\Lambda_{\varepsilon}(\eta)}\bigotimes_{m=1}^d (I_{\alpha_m}^{(m)}-I_{\alpha_m-1}^{(m)}) $ with $I_{\alpha_m}^{(m)}f=\sum_{j=0}^{\alpha_m}f(\chi_j^{(\alpha_m)})\ell_{j}^{(\alpha_m)}(.) $ and $\Lambda_{\varepsilon}(\eta)=\Lambda_{\varepsilon}(\eta_1,\ldots,\eta_d)=\{ (\alpha_1,\ldots,\alpha_d)\in (\mathbb{N}^d)\backslash \{0 \} \ \quad \eta_1^{\alpha_1}\ldots \eta_d^{\alpha_d}>\varepsilon \}.$
\end{remark}


\subsubsection{$L^2$ Convergence rate}


Next, we review error estimates for the stochastic collocation operator in the $L^2$ norm from \cite{AndSchw11,GHL21}.
We start by a Lemma from \cite[p. 1021]{BNT07}.

 \begin{Lemma}\label{lemma_min_pol_approx}
 Let $\Gamma$ be a bounded subset of $\mathbb{R}$ and assume that $f\in \mathcal{C}^{0}(\Gamma,\mathbb{R}) $ admits an analytic extension in the region of the complex plane $ E(\Gamma,\tau):= \{ s\in \mathbb{C} \quad |\quad dist(s,\Gamma)\leq \tau \}$ for some arbitrary $\tau\in \mathbb{R}$.
	 Define $\xi=2\tau/|\Gamma|+\sqrt{1+4\tau^2/|\Gamma|^2}$. Then 
	    \begin{equation}\label{equa5}
	        \min_{p_{r}\in \mathcal{P}_{r}}\|f-p_{r}\|_{L^2(\Gamma)} \leq C\frac{\xi^{-r}}{1-\xi}\|f\|_{L^{\infty}(E(\tau))}
	    \end{equation}
        where C is independent of $r$.
	\end{Lemma}
\begin{Lemma}\label{Lemma_conv_stoch}
Assume that  $f:S\rightarrow \mathbb{R}^n $ admits a complex analytic extension in the region $ E(\tau)=\{ z\in \mathbb{C}^\infty | dist(z_m,[-1;1])< \tau_m \} $ where $ \tau=(\tau_1,\tau_2,\ldots) $ is a sequence of positive numbers such that $ \tau_m\rightarrow \infty $. Define $ \Lambda_\varepsilon(\mu) $ as in \eqref{def_Aeps} where $ \xi_m=\tau_m+\sqrt{1+\tau_m} $ and $ \eta_m=sup_{m^{'}\geq m}\frac{1}{\xi_{m^{'}}} $. Assume that $ \eta_m m^\theta \rightarrow 0 $  for some $ \theta>2(1+\log(4)) $. Consider $H=(\mathbb{R}^n,\|.\|_{\overline{A}})$ where $\|x\|_{\overline{A}}:=\|\overline{A}^{1/2}x\|_{2}$ for all $x\in\mathbb{R}^n$ with $\overline{A}\in\mathbb{S}^{n\times n}_{++}$. Then, for any $ 1>\chi>2(1+\log 4)/\theta $, there exists $ C>0 $ such that 
\begin{equation}
    \|f- \mathcal{I}_{\Lambda_\varepsilon(\eta)}f\|_{L^2(S,H)}\leq C\varepsilon^{1-\chi}\|f\|_{L^\infty(E(\tau),H)}
\end{equation}
for all $ 0<\varepsilon<\eta_1 $.
\end{Lemma}
For the proof of this theorem, we send the reader to  \cite{Bieri09}.
\begin{remark}
    Such convergence rate could also be proved using similar process for any d-variate function. We leave the proof to the reader.
\end{remark}


\section{Error analysis}
\label{Section_approx_av_eig}


Next we study how accurately Ritz eigenvalues from $V$ approximate the eigenvalues $\lambda(\sigma)$ of \eqref{eq:evp} on $(0;\Lambda)$.
By Corollary~\ref{RHS} the accuracy depends on how well eigenvectors related to $(0,\Lambda)$ is approximated in $V$. Recall that $V=\overline{E}_{\leq\rho\Lambda}\oplus X_{\rho\Lambda}$ and any eigenvector on $(0;\Lambda)$ satisfies
\begin{equation}
\label{eq:recall_Z}
    x_k(\sigma)=\overline{x}_k(\sigma) + Z(\sigma,\lambda_k(\sigma))\overline{x}_k(\sigma) \quad \mbox{and} \quad \| x_k(\sigma) \|_M = 1. 
\end{equation}
Then it holds that $\| \overline{x}_k(\sigma) \|_M \leq 1$.
\begin{Def}\label{Def_CL_operator}
    For fixed $\varepsilon>0$, $\eta=(\eta_1,\ldots,\eta_d)\in (0;1)^d$ with $\eta_{j}\geq \eta_{j+1}$ and $q\in \mathbb{N}$, we define the operator $\mathcal{I}_{\Lambda_{\varepsilon}(\eta)}\otimes \pi_q$ applied on $Z$ as
\begin{multline*}
    [(\mathcal{I}_{\Lambda_{\varepsilon}(\eta)}\otimes \pi_q)Z](\sigma,t):= \sum_{\alpha\in \Lambda_{\varepsilon}(\eta)} \sum_{\gamma\leq \alpha}\sum_{j=0}^{q}Z(\chi_{\gamma_1}^{(\alpha_1)},\ldots,\chi_{\gamma_d}^{(\alpha_d)},T_j^{(q)})
     \ell_{1,i}^{(p_i)}(\sigma_1)\cdots\ell_{d,i}^{(p_i)}(\sigma_d)\ell^{(q)}(t)
\end{multline*}
where $\chi_{k}^{(p)}$ is the $k^{th}$ abscissae of the $(p+1)^{th}$ univariate Legendre polynomial and $T_j^{(q)}$ the $k^{th}$ Chebyshev node of degree $q$.
\end{Def}
Corollary~\ref{RHS}  states that 

\begin{equation*}
    \int_{S}\frac{\mu_k(\sigma)-\lambda_k(\sigma)}{\lambda_k(\sigma)} d\sigma
    \leq C \int_{S} \max_{k\in {1,\ldots,N(\sigma)}}  \|(I-P_{A(\sigma)})x_k(\sigma) \|_{A(\sigma)}^2 \; d\sigma,
\end{equation*}
for any $x_k(\sigma)$ such that $\|x_k(\sigma)\|_M=1$. By best approximation property of $P_{A(\sigma)}$ and \eqref{eq:recall_Z}: 
\begin{equation*}
    \|(I-P_{A(\sigma)})x_k(\sigma) \|_{A(\sigma)} 
    = \min_{v \in V} \| \overline{x}_k(\sigma) + Z(\sigma,t)\overline{x}_k (\sigma) - v\|_{A(\sigma)}
\end{equation*}
As $\overline{x}_k(\sigma)\in \overline{E}_{\leq\rho\Lambda}$ 
and $[(\mathcal{I}_{\Lambda_{\varepsilon}(\eta)}\otimes \pi_q) Z]
(\sigma,t)\overline{x}_k(\sigma) \in X_{\rho \Lambda}$, we choose $v = \overline{x}_k(\sigma) +[(\mathcal{I}_{\Lambda_{\varepsilon}(\eta)}\otimes \pi_q )Z]
(\sigma,t)\overline{x}_k(\sigma)$ so that
\begin{equation*}
\|(I-P_{A(\sigma)})x_k(\sigma) \|_{A(\sigma)}^2
\leq
\|[(I-\mathcal{I}_{\Lambda_{\varepsilon}(\eta)}\otimes \pi_q)Z](\sigma,\lambda(\sigma))\overline{x}_k(\sigma)\|_{A(\sigma)}^2.
\end{equation*}
 Then it holds that 
\begin{equation}\label{eq:error_interm}
\int_{S}\frac{\mu_k(\sigma)-\lambda_k(\sigma)}{\lambda_k(\sigma)} d\sigma
\leq C \int_{S} \max_{\substack{y\in\overline{E}_{\leq\rho\Lambda}\\ \|y\|_{M}\leq 1} } \|[(I-\mathcal{I}_{\Lambda_{\varepsilon}(\eta)}\otimes \pi_q)Z](\sigma,\lambda(\sigma)) y \|_{A(\sigma)}^2 \; d\sigma.
\end{equation}

Hence, upper bound for Ritz eigenvalue error follows by estimating the interpolation error in \eqref{eq:error_interm}. 


\subsection{Convergence rate}


In the rest of the paper for fixed $\sigma\in S$, we denote $Z_{\sigma}(t):=Z(\sigma,t)$ for all $t\in (0,\Lambda)$ and for any fixed $t\in (0;\Lambda)$, $Z_t(\sigma)=Z(\sigma,t)$ for all $\sigma\in S$.
Also, let $\mathcal{H}(\sigma)$ be the Hilbert space $\mathbb{R}^n$ equipped with $A(\sigma)$  inner product and the induced norm. Let $F:(0;\Lambda)\mapsto \mathcal{H}(\sigma)$. Then we use the notation 
\begin{equation}
\|F \|_{L^\infty\left( (0;\Lambda) ; \mathcal{H}(\sigma) \right) } = \sup_{t\in (0;\Lambda)} \| F(t) \|_{A(\sigma)}
\end{equation}
With these notations, we now proceed to give technical results needed to bound the RHS of \eqref{eq:error_interm}.   
\label{subsect_CV_rate}

For fixed degree $q$ and $t\in (0;\Lambda)$, denote the operator $Z_{q,t}(\sigma)=\pi_q Z(\sigma,t)$ for all $\sigma\in S$ where $\pi_q$ is the Lagrange interpolation operator in the t-variable at Chebyshev points of degree $q$. 

Note that $(\mathcal{I}_{\Lambda_{\varepsilon}(\eta)}\otimes \pi_q) Z = \mathcal{I}_{\Lambda_{\varepsilon}(\eta)} Z_{q,t}$. One can easily see using separation in the interpolation that for all $y\in\overline{E}_{\leq \rho\Lambda}$,
\begin{multline}\label{equ:sep_norm}
    \frac{1}{2}\|Z(\sigma,t)y-[(\mathcal{I}_{\Lambda_{\varepsilon}(\eta)}\otimes \pi_q)Z](\sigma,t)y\|_{A(\sigma)}^2
    \leq  \|[(I-\pi_q) Z_{\sigma}](t)y\|_{A(\sigma) } ^2\\
    + \|[(I-\mathcal{I}_{\Lambda_{\varepsilon}(\eta)}) Z_{q,t}](\sigma)y\|_{A(\sigma)}^2
\end{multline}
    and we then study separately these quantities to obtain our main result in theorem \ref{Thm_Bauer_Fikes_like}.
As one can observe that $\|[(I-\pi_q)Z_{\sigma}](t)y\|_{A(\sigma)}\leq \|(I-\pi_q) Z_\sigma y\|_{L^\infty\left( (0,\Lambda); \mathcal{H}(\sigma) \right) } $, we first obtain:
 \begin{Lemma}\label{Lemma_first_bound_pi_q}
     Let $\sigma\in S\subset \mathbb{R}^d$, and $A:S\rightarrow \mathbb{S}^{n\times n}_{++}$. Assume that $A$ is spectrally equivalent to $\overline{A}\in\mathbb{S}^{n\times n}_{++}$ with constants $\alpha$ and $\beta$ as in \eqref{eq:spequiv}. In addition, let $\rho>2$ and $\pi_q$ be the Lagrange interpolant at Chebyshev points of order $q\in\mathbb{N}$. 
     Then it holds that
     \begin{equation}\label{upper_bnd_pi_q_approx}
          \|(Z_\sigma-\pi_q Z_\sigma )y\|_{L^\infty\left( (0,\Lambda); \mathcal{H}(\sigma) \right) }^2
         \leq \frac{\beta\rho^3}{4^{2q-1}(\rho-1)^{2(q+1)}}\Lambda
     \end{equation}
     for any $y\in \overline{E}_{\leq \rho\Lambda}$ such that $\|y\|_{M}\leq 1$. 
 \end{Lemma}
 \begin{proof} The proof utilises~Lemma \ref{Lemma_express_Z} to relate the interpolation error of $Z(\sigma,t)$ to interpolation error of a scalar function. This diagonalisation argument allows us to obtain estimates independent on the dimension of the matrix $Z(\sigma,t)$. Let $y\in \overline{E}_{\leq\rho \Lambda}$ and $(\xi_k,u_k)$, $u_i^T M u_j = \delta_{ij}$ be the eigenmodes of the pencil $(A_\perp(\sigma),M_\perp)$. By Lemma \ref{Lemma_express_Z}
 \begin{equation*}
     Z_\sigma(t)y = W_{\perp}\sum_{k} u_k(\sigma) f_k(t) \beta_k \hspace{0.2cm} \mbox{for} \hspace{0.2cm}  \beta_k = u_k(\sigma)^{T}W_{\perp}^{T}A(\sigma)y \hspace{0.2cm} \mbox{and} \hspace{0.2cm} f_k(t) = (\xi_k(\sigma)-t)^{-1}.
 \end{equation*}
 Hence, $\pi_q Z_{\sigma}(t) = \sum_{k} u_k(\sigma) \pi_q f_k(t) \beta_k$. By orthogonality properties of the eigenbasis $(u_k(\sigma))_k$ it holds that
\begin{equation}
\label{eq:Zy_norm}
        \| Z_{\sigma}(t) y\|^2_{A(\sigma)} = \sum_{k}  \xi_k(\sigma) f_k(t)^2 \beta^2_k 
        \end{equation}

and
    \begin{equation}
    \label{eq:Zy_iperror_norm}
    \| Z_{\sigma}(t)y-(\pi_q Z_{\sigma})(t)y\|^2_{A(\sigma)}=\sum_k \xi_k(\sigma) (f_k(t)-\pi_q f_k(t))^2 \beta_k^2.
    \end{equation}
    As $f_k$ is monotonically increasing and  by \eqref{eq:Zy_norm} it follows that $$\|Z_\sigma y\|_{L^\infty\left( (0,\Lambda); \mathcal{H}(\sigma) \right) }^2 = \sum_k \xi_k(\sigma) f_k(\Lambda)^2 \beta^2_k.$$ 
    We proceed to connect this expression to \eqref{eq:Zy_iperror_norm}. By well-known Lagrange-Chebychev interpolation error estimate it holds that 
    \begin{equation}
        \max_{t\in(0;\Lambda)}|f_k - \pi_q f_k|\leq \frac{\Lambda^q}{2^{2q-1} (q+1)!}\max_{t\in (0,\Lambda)}|f^{(q+1)}_k(t)|=\frac{2\Lambda^q}{4^q } f_k(\Lambda)^{q+1},
    \end{equation}
    see \cite{Stewart96}. Hence,
    \begin{equation}
    \| Z_{\sigma}(t)y-(\pi_q Z_{\sigma})(t)y\|^2_{A(\sigma)} \leq \sum_k \frac{4\Lambda^{2q}}{4^{2q}} f_k(\Lambda)^{2q} \xi_k(\sigma) f_k(\Lambda)^2  \beta_k^2.
    \end{equation}

    As $\xi_k(\sigma) \geq \rho \Lambda$ and $t\in(0,\Lambda)$, one can easily check that $\Lambda f_k(\Lambda) \leq (\rho-1)^{-1}$. Hence,
    \begin{equation}
    \| Z_{\sigma}(t)y-(\pi_q Z_{\sigma})(t)y\|^2_{A(\sigma)} \leq \frac{1}{4^{2q-1}(\rho-1)^{2q}} \sum_k \xi_k(\sigma)  f_k(t)^2  \beta_k^2.
    \end{equation}
    Application of \eqref{eq:Zy_norm} and Lemma \ref{Lemma_bound_Zy} completes the proof. 
    \end{proof}

    \begin{Lemma}\label{Lemma_second_bound_I_Lambda}
        Let $t\in (0,\Lambda)$. Under the conditions of Lemma \ref{Lemma_first_bound_pi_q}, for $\varepsilon>0$, $1>\eta_1\geq \eta_2\geq \ldots \geq \eta_d >0 $ with $d\in \mathbb{N}$, there exists $C=C(\alpha,\beta)$ such that
    
    \begin{equation}
        \| ( I - \mathcal{I}_{\Lambda_{\varepsilon}(\eta)} )Z_{q,t} y \|_{L^2(S,\mathcal{H}(\sigma))}^2
     \leq C \varepsilon^{2(1-\mathcal{X})}\Big(1+\frac{1}{4^{2q-1}(\rho-1)^{2q}}\Big)\Big(\frac{\rho}{\rho-1}\Big)^2\rho\Lambda 
    \end{equation}
    for any $y\in \overline{E}_{\leq \rho\Lambda}$ with $\|y\|_{M}\leq 1$ and $ 1>\mathcal{X}>2(1+\log 4)/\theta $ when $\eta_m\leq (1+m)^{-\theta}$ with $\theta>1$.
    \end{Lemma}
    \begin{proof}
        First observe that the spectral equivalence implies, $$\|( Z_{q,t}-\mathcal{I}_{\Lambda_{\varepsilon}(\eta)} Z_{q,t})y\|^2_{A(\sigma)}
    \leq \beta \|(\overline{A}^{1/2} Z_{q,t}-\mathcal{I}_{\Lambda_{\varepsilon}(\eta)} \overline{A}^{1/2}Z_{q,t})y\|^2_{2}$$
    for any $y\in \overline{E}_{\leq \rho\Lambda}.$
    From the $L^2(S)$-approximation error estimate for the stochastic collocation operator given in Lemma \ref{Lemma_conv_stoch} it follows that 
    \begin{equation}\label{stoch_uppbd_Abar}
        \|(\overline{A}^{1/2}Z_{q,t}-\overline{A}^{1/2}\mathcal{I}_{\Lambda_{\varepsilon}(\eta)}Z_{q,t})y\|_{L^2(S,(\mathbb{R}^n,\|.\|_2))}
        \leq \varepsilon^{2(1-\mathcal{X})} \max_{\sigma \in  S} \| \overline{A}^{1/2} Z_{q,t}(\sigma) y\|_2^2.
    \end{equation}
    By spectral equivalence
    \begin{equation*}
       \max_{\sigma \in  S} \| \overline{A}^{1/2} Z_{q,t}(\sigma) y\|_2^2   \leq\frac{1}{\alpha}\max_{\sigma\in S} \|Z_{q,t}(\sigma) y\|_{A(\sigma)}^2
    \end{equation*} 
    for any $y\in \overline{E}_{\leq\rho\Lambda}$. Then, as $\|Z(\sigma,t)y\|_{A(\sigma)}\leq \frac{\rho}{\rho-1}(\beta\rho\Lambda)^{1/2}\|y\|_{M}$ by Lemma~\ref{Lemma_bound_Zy} and using decomposition $ \|Z_{q,t}(\sigma) y\|_{A(\sigma)}^2\leq\|\pi_q Z_{\sigma}y\|_{L^{\infty}((0;\Lambda),\mathcal{H}(\sigma))}^2\leq 2\|Z_{\sigma}y\|_{L^{\infty}((0;\Lambda),\mathcal{H}(\sigma))}^2+2\|(I-\pi_q)Z_{\sigma}y\|_{L^{\infty}((0;\Lambda),\mathcal{H}(\sigma))}^2 $ for all $\sigma \in S$ combined with \eqref{upper_bnd_pi_q_approx}, \eqref{stoch_uppbd_Abar} and Lemma \ref{Lemma_bound_Zy}, we finally obtain the claimed result.
    \end{proof}
%

\begin{Thm}\label{Thm_Bauer_Fikes_like}
    Under the conditions of Lemmas \ref{Lemma_first_bound_pi_q} and \ref{Lemma_second_bound_I_Lambda}, denote by $\lambda_k(\sigma)$ the exact $k^{th}$ eigenvalue of \eqref{eq:evp} associated with parameter $\sigma\in S$ and $\mu_k(\sigma)$ their approximation from the Ritz space $V$. If for all $\sigma\in S$,  $e(\sigma,V)<1/2$ with $e(\sigma,V)$ defined in Corollary \ref{RHS}, then for any $\varepsilon>0$ and for any $ 1>\chi>2(1+\log 4)/\theta $,
    
     \begin{equation}\label{Bauer_Fike_upper_bound}
        \int_S \max_{k\in \{1,\ldots,N(\sigma)\} } \frac{\mu_k(\sigma)-\lambda_k(\sigma)}{\lambda_k(\sigma)} \; d\sigma 
    \leq 
    C\Big(\frac{\rho^2}{4^{2q-1} (\rho-1)^{2q}}+2\varepsilon^{2(1-\chi)}\frac{\beta}{\alpha}\Big)\frac{\rho^3}{(\rho-1)^2}\Lambda
    \end{equation}
    with $C=C(\beta,\alpha)$.
    \end{Thm}
\begin{proof}
    Recall from equations \eqref{eq:error_interm} and \eqref{equ:sep_norm} that 
    \begin{multline*}
        \int_S \max_{k\in \{1,\ldots,N(\sigma)\} } \frac{\mu_k(\sigma)-\lambda_k(\sigma)}{\lambda_k(\sigma)} \; d\sigma 
        \leq 2  \int_{S}\max_{\substack{y\in\overline{E}_{\leq\rho\Lambda}\\ \|y\|_{M}\leq 1}}\|(Z_\sigma-\pi_q Z_\sigma )y\|_{L^\infty\left( (0,\Lambda); \mathcal{H}(\sigma) \right) }^2 d\sigma\\
    + \max_{\substack{y\in\overline{E}_{\leq\rho\Lambda}\\ \|y\|_{M}\leq 1}}\|(I-\mathcal{I}_{\Lambda_{\varepsilon}(\mu)}) Z_{q,t} y \|_{L^2(S;\mathcal{B}(\mathcal{H}(\sigma))}^2.
    \end{multline*}
    Thus
    \begin{multline*}
        \int_S \max_{k\in \{1,\ldots,N(\sigma)\} } \frac{\mu_k(\sigma)-\lambda_k(\sigma)}{\lambda_k(\sigma)} \; d\sigma 
        \leq C\Big(  |S| \frac{\rho^2}{4^{2q-1}(\rho-1)^{2(q+1)}} \\
        + \varepsilon^{2(1-\mathcal{X})}\frac{\beta}{\alpha}\Big(1+\frac{\rho}{4^{2q-1}(\rho-1)^{2q}}\Big)\frac{\rho^3\Lambda}{(\rho-1)^2}\Big)
    \end{multline*}
    from decomposition and upper bounds in Lemmas \ref{Lemma_first_bound_pi_q} and \ref{Lemma_second_bound_I_Lambda}, we obtain the claimed result.
\end{proof}

\section{Numerical experiments}


\label{num_exp_section}

In this section, we demonstrate the proposed Ritz method by solving GEVP arising from finite element discretization of the problem: for given $(a_i)_{i=0}^d \subset L^\infty(\Omega)$, $\sigma \in S$ find $(\lambda,u) \in \mathbb{R}_{+} \times H^1_0(\Omega) \setminus \{0 \}$ s.t. 
\begin{equation}
\label{eq:example}
    \sum_{i=0}^d \sigma_i 
 \int_{\Omega} a_i  \nabla u \cdot \nabla v \; dx = \lambda \int_{\Omega} u v \; dx \quad \mbox{for any $v\in H^1_0(\Omega)$}. 
\end{equation} 
This problem is similar to the one obtained from stochastic PDE after K-L expansion, see  Example~\ref{exemple2.2}. After finite element discretization, \eqref{eq:example} with properly chosen coefficient functions is of form \eqref{eq:evp} for some $A_0,M \in \mathbb{S}^{n\times n}_{++}$ and $A_1,\ldots,A_d \in \mathbb{S}^{n\times n}$. In this section, the discretization is done using an in-house finite element code that uses standard first order triangular Lagrange elements, see \cite{Braess07}.  All examples are computed using the same mesh that is obtained by uniform refinement and leads to $n=3969$.

We propose two different algorithms for computing a basis matrix for $V$. First one simply computes a basis matrix from the Chebyshev-Legendre interpolation operator  and the second one uses additional low-rank approximation and deflation steps. 

Let $\widehat{S} \subset S$ be a finite sample set. In the following examples, we numerically estimate the approximation error of our method as
\begin{equation}\label{Error_studied}
    \max_{\sigma \in \widehat{S} } \max_{ i \in \{1,\ldots, N  \}} \frac{\mu_i(\sigma)-\lambda_i(\sigma)}{\lambda_i(\sigma)}
\end{equation}
where $\lambda_i(\sigma)$ is the $i^{th}$ smallest eigenvalue of \eqref{eq:evp} at $\sigma$ and $\mu_i(\sigma)$ its Ritz approximation.  For simplicity, we use fixed number of eigenvalues $N$ instead of interval $(0;\Lambda)$.


\subsection{Ritz basis construction}

%
In this section, we give two methods for computing the basis matrix of the Ritz space $V=\overline{E}_{\leq\rho\Lambda}\oplus X_{\rho\Lambda}$.
Let $\{\chi_i^{(\alpha)}\}\subseteq S$, $\{T_j^{(q)}\}\subseteq (0;\Lambda)$ be collocation points. Recall that $W=[w_1,\ldots,w_m]$ is the basis matrix of $\overline{E}_{\leq\rho\Lambda}$, the sample vectors are defined as $s_{ijk}= Z(\chi_i^{(\alpha)},T_j^{(q)})w_k$, and the subspace $X_{\rho\Lambda}$ is given by $X_{\rho\Lambda} :=span_{i,j,k} s_{ijk}$. A simple way to obtain a basis matrix $Q$ of $V$ is to collect basis vectors of $W$ and all sample vectors as columns of an auxiliary matrix and apply the $qr$-decomposition to obtain unitary basis matrix $Q$ of $V$. 

This method, which we refer to as RMCLI (\textit{Ritz Method with Chebyshev-Legendre Interpolation}), has two challenges. First, we collect basis vectors of $W$ and all sample vectors as columns of an auxiliary matrix. For a large number of interpolation points, this matrix may become so large that issues in terms of storage arise and computing its $qr$-decomposition becomes too costly. Second, as our aim is to reduce the computational cost of approximating eigenvalues, the size of subspace $V$ should be as small as possible. 

 To avoid these issues, we propose a modified version of the previous practical implementation to reduce the memory requirement at each step, avoid the $qr$-decomposition, and further reduce the dimension of $V$. We refer to this new implementation as RMCLI reduced. In the reduced version, we construct a sequence subspaces $V^{(k)}$ as follows. The first space $V^{(0)} = \overline{E}_{\leq\rho\Lambda}$ with orthonormal basis $Q_0$. We then evaluate SVD of $(I-Q_{k-1}Q_{k-1}^T)Z(\chi_i^{(\alpha)},T_j^{(q)})W$, that represents additional information from the current collocation point to the Ritz space. The right singular vectors corresponding to singular values larger than given tolerance are then added to the space $V^{(k-1)}$ to obtain the space $V^{(k)}$. The final space obtained by repeating this process for all interpolation points is the Ritz space of RMCLI reduced. This procedure avoids storing large auxiliary matrices. In addition, the low-rank approximation reduces the dimension of the produced Ritz space. However, we do not provide analysis on the effect of the low-rank approximation to the eigenvalue error. 


To use either of the two RMCLI variants one first needs to fix values of the oversampling parameter $\rho$, Chebychev interpolation degree $q$, and stochastic collocation parameters $\varepsilon,\eta$. In our numerical tests we found that choosing $\rho \in (1;1.5)$ yields Ritz spaces with small rank and at the same time gives sufficiently accurate results. Moreover, as it will be shown in example \ref{example_epsilon}, the choice $\varepsilon := \eta_1/2$ and small values of $q$ allowed us to reach sufficient accuracy while keeping the computational cost small.


\subsection{One dimensional parameter} 


   In this section, we consider the problem \eqref{eq:example} with $a_0(x)=1$ and $a_1(x)=\frac{1}{2}\sin(\pi x)$. We start by studying the impact of the parameter $\varepsilon$ on the number of interpolation points for $\sigma$-parameter and the associated Ritz eigenvalue error.
\begin{example}\label{example_epsilon}
\normalfont In the one dimensional case, the stochastic collocation part of the operator $\mathcal{I}_{\Lambda_{\varepsilon}(\eta)}\otimes \pi_q$  on $ S\times (0;\Lambda)=[-1;1]\times (0;\Lambda) $  consists in interpolating at $\alpha_{\max}+1=\max\{\alpha| \alpha\in\Lambda_{\varepsilon}(\eta)\}+1$ points. Indeed, when $d=1$, formula \eqref{def_inf_Stoch_op} is reduced to the Lagrange polynomial $I_{\alpha_{\max}}$ which uses points $\chi_0^{(\alpha_{\max})},\ldots,\chi_{\alpha_{\max}}^{(\alpha_{\max})}$. In Table \ref{tab:eps_chi}, the number of interpolation points used on $\sigma$-variable using RMCLI and maximum value of \eqref{Error_studied} for the 10 smallest eigenvalues with $\eta=0.5$, $q=2$, and different values of $\varepsilon<\eta$ is displayed. The error was evaluated using equidistant sample set with $100$ points for $\sigma$. 

We first observe that values of $\varepsilon$ close to $\eta$ already provide sufficient accuracy for many applications. In particular, choosing $\varepsilon = \eta/1.1$ provides a $10^{-5}$ relative error and has only two interpolation points $\chi_{i}^{(\alpha)}\in S$. Moreover, as we reduce $\varepsilon$, the number of interpolation points in $\sigma$ stays small while accuracy improves quickly with error reaching $10^{-10}$ precision using only six $\sigma$-interpolation points. In this example our method gives highly accurate results with a very small number of interpolation points, and thus, evaluations of the correction operator. 

\begin{table}[h]
        \centering
        \begin{tabular}{|c|c|c|c|c|}
        \hline
            $\varepsilon$ & $\eta/1.1$ & $\eta/2$ & $\eta/5$ & $\eta/20$ \\
             \hline
            Number of stochastic points $\chi_{i}^{(\alpha)}$ & 2 & 3 & 4 & 6 \\
             \hline
             maximum value of \eqref{Error_studied} & $9,664. 10^{-6}$ & $7,465. 10^{-6}$ & $1,254. 10^{-8}$ & $1,208.10^{-10}$ \\
             \hline
        \end{tabular}
        \caption{Number of $\sigma$-interpolation points in the stochastic collocation interpolation for different values of $\varepsilon$.}
        \label{tab:eps_chi}
    \end{table}
    
    We fix $\eta=0.5$ and $\varepsilon=\eta/10$, impose a tolerance of $10^{-7}$ in RMCLI reduced and measure \eqref{Error_studied} for degree $q$ from 1 to 5 for the Chebyshev interpolation in the eigenvalue parameter. In Table \ref{table_dim_V}, we measure the dimension of the Ritz spaces $V$ produced by RMCLI and RMCLI reduced for different number of Chebyshev interpolation points $q$. We see the efficiency of RMCLI reduced in reducing the dimension of subspace obtained and speeding up the solution of the reduced eigenvalue problem as the dimension of the Ritz subspace grows from $60$ to $260$ as $q$ increases from $1$ to $5$. 

\begin{table}[h]
     \centering
     \begin{tabular}{|c|c|c|c|c|c|c|}
     \hline
    $q$ & 1  & 2 & 3 & 4 & 5 \\
     \hline
     $\mathop{dim} (V) $ w/ RMCLI & $60$ & $110$ & $160$ & $210$ & $260$ \\
     \hline
     $\mathop{dim} (V) $ w/ RMCLI reduced & $50$ & $90$ & $120$ & $139$ & $151$ \\
        \hline 
     \end{tabular}
     \caption{Dimensions of Ritz subspaces generated by RMCLI and RMCLI reduced for $\eta=0.5$, $\varepsilon=\eta/10$ and tolerance equal to $10^{-7}$.}
     \label{table_dim_V}
 \end{table}

 Moreover, in Figure \ref{fig:red_basis}, we see that even with a $10^{-7}$ tolerance in RMCLI reduced, we approximately reach the same accuracy as with RMCLI. Thus, we would advise to use RMCLI reduced as it reaches similar accuracy while keeping the dimension of the Ritz subspace small.

 \begin{figure}[H]
    \centering
    \includegraphics[scale=0.6]{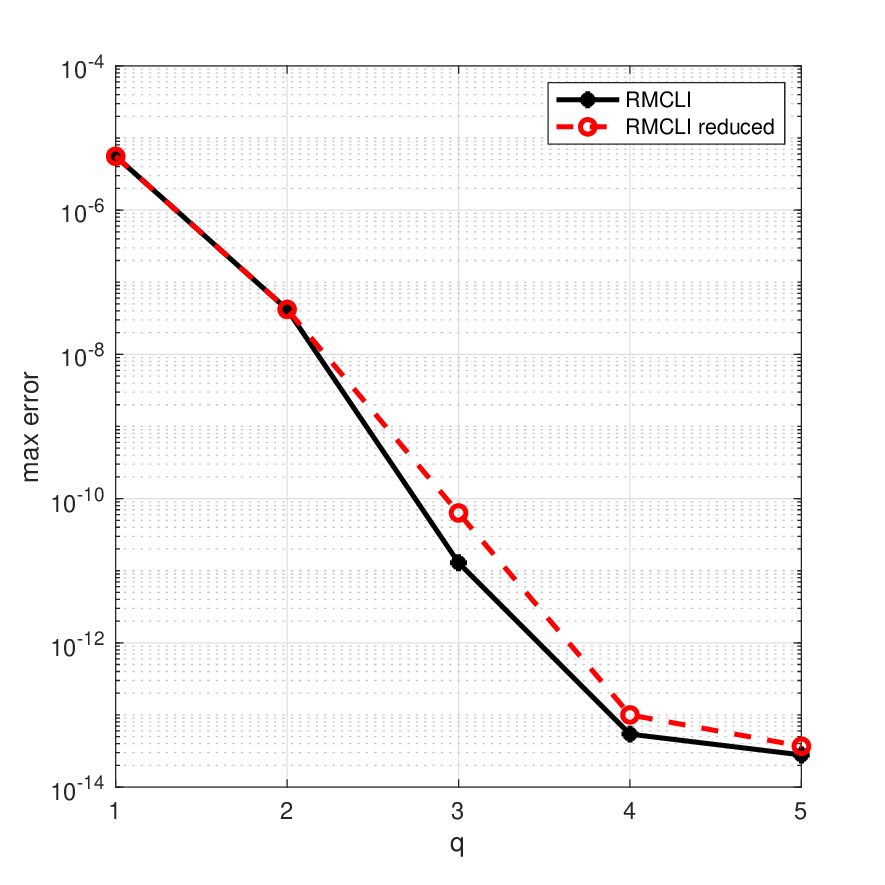}
    \caption{Values of \eqref{Error_studied} for the $10$ smallest exact eigenvalue solutions using RMCLI (black line),  and reduced version (red line) with $\eta=0.5 $,  $\varepsilon=\eta/10$ and $q=1,\ldots,5$.}
    \label{fig:red_basis}
    \end{figure}
\end{example}

\begin{example}  \normalfont We continue by comparing the Chebyshev-Legendre points used to construct our Ritz subspace to other interpolation points in the plane. Ritz eigenvalue error \eqref{Error_studied} with different interpolation points is given Figure \ref{fig:comp_interp_pts}. The basis matrix was constructed using RMCLI while replacing the Chebyshev-Legendre points with others. Three types of points were considered:
\begin{enumerate}
\item Multivariate Legendre interpolation points $ (s_i,t_j)=(\chi_{i}^{(\alpha_1)},\chi_{j}^{(\alpha_2)})\in S\times [0;\Lambda]$ where $\chi_{\ell}^{\alpha_{k}} $ is the $\ell^{th}$ root of the $\alpha_k^{th}$ Legendre polynomial for $\ell\in \{ 0,\ldots,\alpha_k \}$, $k=1,2$. Note that this interpolation is equivalent to use the Stochastic collocation operator described in \ref{subsubsection_stoch_colloc_def} for both $\sigma$ and $t$-variables.
    \item Chebyshev-Legendre points $(s_i,t_j)=(\chi_{i}^{(\alpha)},T_j^{(q)})\in  S\times [0;\Lambda]$ with $\chi_i^{(\alpha)}\in S$ as before and $ T_j^{(q)}=  \frac{\rho\Lambda}{2} + \frac{\rho\Lambda}{2}\cos \big( \frac{2j+1}{2q}\pi \big{)} $ the Chebyshev points of order $q$. This leads to the interpolation operator $\mathcal{I}_{\Lambda_{\varepsilon}(\mu)}\otimes \pi_q$.
    \item Padua points $(s_i,t_j)=(s_i^n,t_j^n)\in S\times [0;\Lambda]$, defined as the equally spaced points of the parametric curve of degree $n$ $ \gamma_n(t)=(\cos(nt),\cos((n+1)t)) $\cite{CDMV08}. If we denote $s_i^n=\frac{\sigma_{min}+\sigma_{max}}{2} + \frac{\sigma_{max}-\sigma_{min}}{2}\cos(\frac{i\pi}{n})$ and $t_j^n=\frac{\Lambda}{2} + \frac{\Lambda}{2}\cos(\frac{j\pi}{n})$, then the Padua points are given explicitly  by 
    \begin{itemize}
        \item if $n=2m$, $s_i = (s_{2i_1+1}^n,s_{2i_2}^{n+1})$ for $ 0\leq i_1\leq m-1, 0\leq i_2\leq m $ and $t=(t_{2j_1}^n,t_{2j_2 +1}^{n+1})$ for $ 0\leq j_1\leq m, 0\leq j_2\leq m $;
        \item if $n=2m+1$, $s_i = (s_{2i_1+1}^n,s_{2i_2}^{n+1})$ for $ 0\leq i_1\leq m, 0\leq i_2\leq m+1 $ and $t=(t_{2j_1}^n,t_{2j_2 +1}^{n+1})$ for $ 0\leq j_1\leq m, 0\leq j_2\leq m $;
    \end{itemize}
    \end{enumerate}
    \begin{figure}[H]
    \centering
    \includegraphics[scale=0.6]{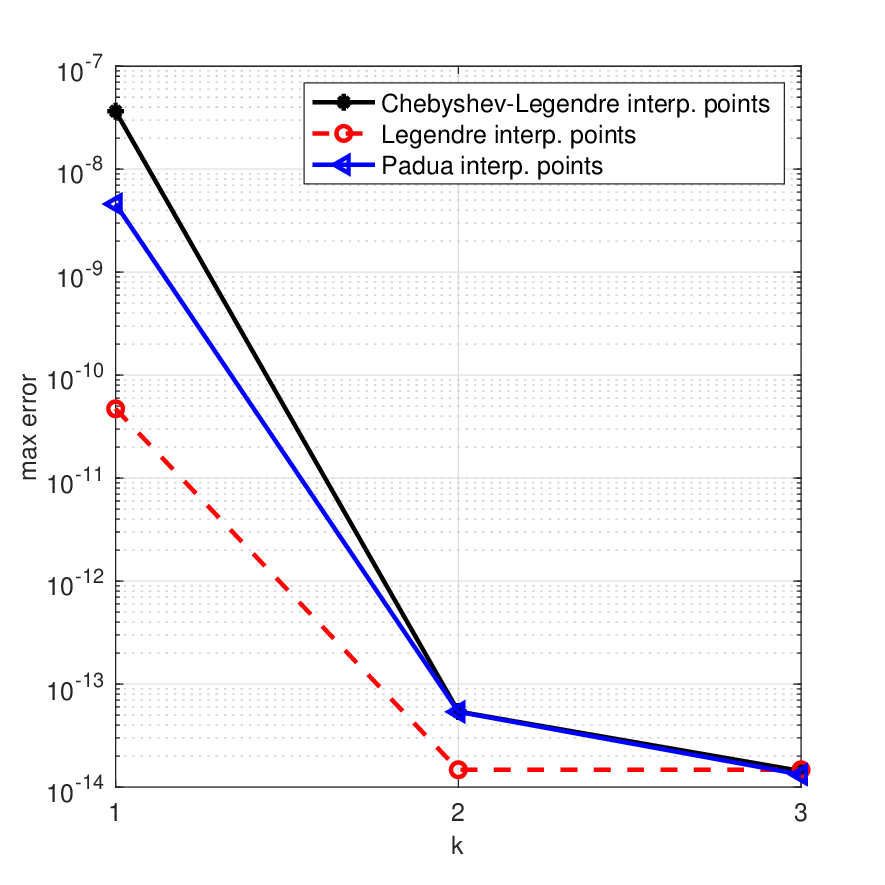}
    \caption{Value of \eqref{Error_studied} for the $10$ smallest exact eigenvalues of \eqref{eq:evp} for $10$ values of $\sigma$ using RMCLI with $\eta=0.5$,  $\varepsilon = \eta/(10 k)$ and $q=2k$ (black line), Legendre points with $\eta=(0.5,0.25)^T$ and $\varepsilon = \eta_1/(10 k)$ (red dashed line) and Padua points (blue line) for $\rho\Lambda=320$, after $6$ refinements of mesh for $k=1,2,3$.}
    \label{fig:comp_interp_pts}
\end{figure}
In Figure \ref{fig:comp_interp_pts}, we measure Ritz eigenvalue errors \eqref{Error_studied} obtained from CL interpolation with the method subspace constructed using three different types of interpolation points $(s_i,t_j)$. 
In the case of Chebyshev-Legendre interpolation points, we measure this error for $q=2k$ with $k=1,2,3$ the Chebyshev degree in the last parameter interpolation, $\mu=0.5$ and $\varepsilon =\eta/(10k)
$. For the Multivariate Legendre case, we consider $\mu=(0.5 \hspace{0.2cm}0.25)^T$ and then $\varepsilon =\eta_1/(10k) $ for $k=1,2,3$ and for Padua points we consider $n=3\times k$.

Based on Figure~\ref{fig:comp_interp_pts},  Legendre, Chebyshev-Legendre and Padua points all yield equivalent accuracy. Namely, for $k=3$ we reach the same accuracy of $10^{-14}$ using RMCLI with Chebyshev-Legendre, multivariate Legendre points or Padua points. 

Note that in Table \ref{tab:my_label}, the number of Padua interpolation points used for $k=1,2$ and $3$ are $10$, $28$ and $55$, respectively, while we use $10$, $20$ and $30$ Chebyshev-Legendre interpolation points consecutively and $19$, $28$ and $28$ Legendre interpolation points. Also, while using the Legendre interpolation points, the dimension of Ritz subspace stays unchanged for $k=2$ and $k=3$. Indeed, multi-index set $\Lambda_{\varepsilon}(\eta)$ is the same for both cases as the value of $\varepsilon$ is not decreased considerably to have more indices $\alpha$ in $\Lambda_{\varepsilon}(\eta)$ and thus the number of interpolation points is the same.

\begin{table}[h]
    \centering
    \begin{tabular}{|c|c|c|c|c|c|c|c|}
    \hline
       k  & 1 & 2 & 3 & & 1 & 2 & 3 \\
         \hline
       num. of. sample vectors, CL  & 110 & 210 & 310 & $\dim(V)$, CL  &  90 & 170 & 202  \\
         \hline
        num. of. sample vectors, L  & 200 & 290 & 290 & $\dim(V)$, L  & 110 & 182 & 182 \\
         \hline
       num. of. sample vectors, Padua  &  110 & 290 & 560 & $\dim(V)$, Padua  & 110 & 200 & 225 \\
         \hline
    \end{tabular}
    \caption{Dimensions of Ritz subspaces before and after orthogonalisation generated by RMCLI.}
    \label{tab:my_label}
\end{table}
In this example, different interpolation points yield similar results. Namely, the dimensions of the different Ritz spaces are almost identical. However, we have only proved convergence for the CL-family of the points that are used also in the multidimensional cases.

\end{example}


\subsection{The multidimensional case}


We now turn to the multidimensional case and consider the problem \eqref{eq:example} with $a_0(x)=1$, $a_m(x)=2^{-m}\sin (m\pi x)$ for $m\geqslant 1$, $S=[-1;1]^d$, $d\in\{1,2,3,4\}$. The index set for stochastic collocation $\Lambda_{\varepsilon}(\eta_1,\ldots,\eta_d)$ is chosen such that $\eta_m\leq(1+m)^{-\theta}$ for some $\theta>1$. 
For numerical reasons, following \cite{GHL21}, we replace $\Lambda_{\varepsilon}(\eta)$ with $\tilde{\Lambda}_{\varepsilon}(\eta)=\{ \alpha\in\Lambda_{\varepsilon}(\eta) : c_{\alpha}\neq 0 \}$ with $c_{\alpha}=\sum_{\gamma\in\Lambda_{\varepsilon}(\eta)}\mathds{1}_{\{\gamma-1\leq \alpha\leq \gamma\}} (-1)^{\|\alpha-\gamma\|_{\ell^1(\mathbb{N})}}$ that avoid computing unnecessary terms.

\begin{example}\label{example_q}
\normalfont
We first study the Ritz eigenvalue approximation error with respect to parameter dimension $d$ and Chebyshev degree $q$. Figure \ref{fig:comp_dim_and_q} displays the error measure \eqref{Error_studied} when the Ritz eigenvalues are computed using RMCLI and its reduced version with a tolerance of $10^{-6}$ for $10$ equidistant values of $\sigma\in S$ with Chebyshev-Legendre interpolation points and for Chebyshev interpolation  degree $q=1,3$ and $5$. We observe that regardless of the SVD reduction step in RMCLI, both methods reach similar precision.
Note that in all cases, the error grows with dimension and increasing the number of Chebyshev interpolation points does not considerably affect the error. This is the case, especially in dimension $d=4$. The authors believe that in this case, the Stochastic collocation part of the interpolation is restricting the accuracy. This is, $\eta_m$ should be decreased to reach more accurate results. 

\begin{figure}[H]
    \centering
    \includegraphics[scale=0.5]{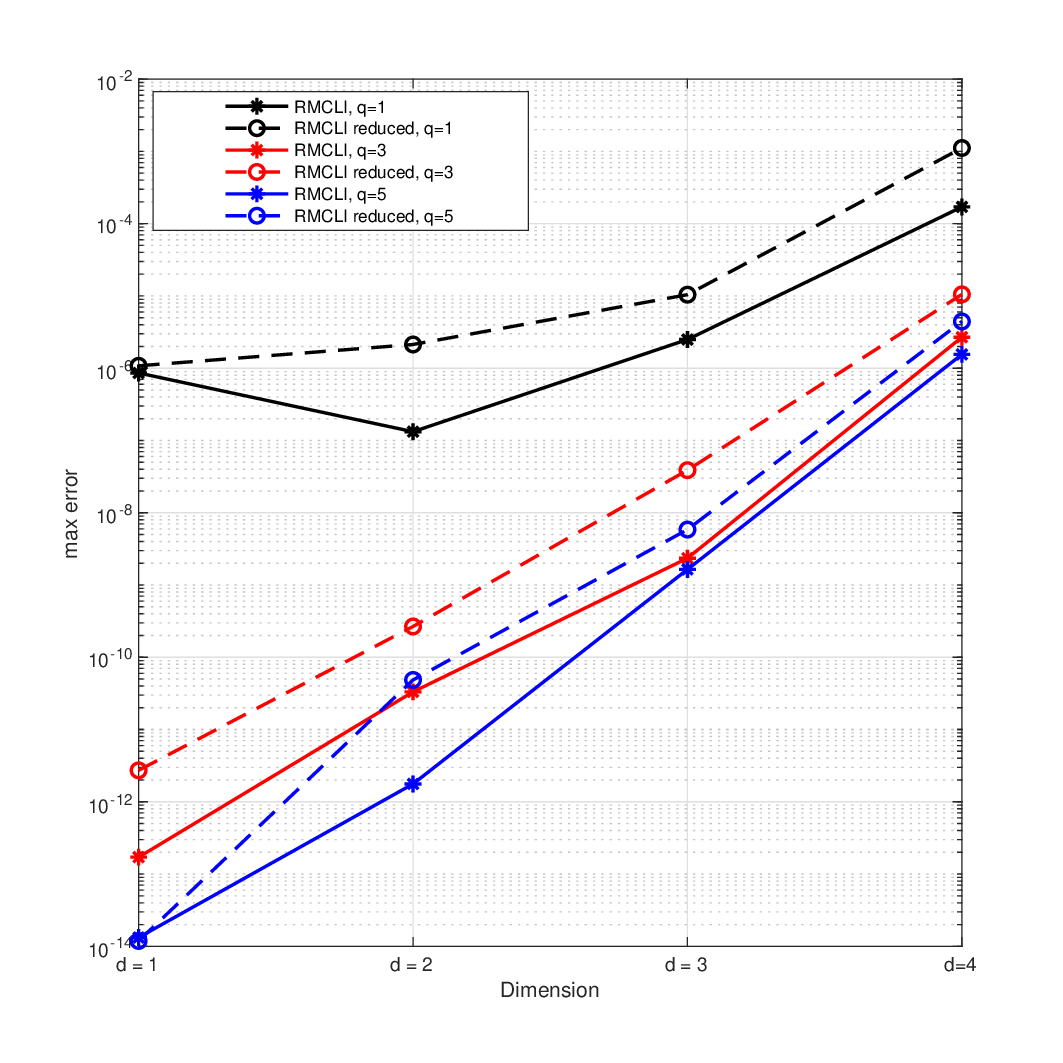}
    \caption{Value of \eqref{Error_studied} in the $1,2,3$ and $4$d cases with $S=[-1;1]^d$ for the $10$ smallest eigenvalues with $\lambda_i(\sigma)$ for 10 value of $\sigma$ using RMCLI and its reduced version for $q=1,3$ and $5$.}
    \label{fig:comp_dim_and_q}
\end{figure}

\begin{table}[h]
    \centering
    \begin{tabular}{|c|c|c|c|c|}
    \hline
        parameter dimension d & 1 & 2 & 3 & 4 \\
        \hline 
        $\dim (V)$ w/ RMCLI, $q=1$ & 101 & 149 & 192 & 209\\
        \hline
         $\dim (V)$ w/ RMCLI, $q=3$ & 197 &  280 &  317  & 324\\
        \hline
        $\dim (V)$ w/ RMCLI, $q=5$ & 252 &  322  & 305 &  314\\
        \hline
        $\dim (V)$ w/ RMCLI reduced, $q=1$ & 79  & 124  & 152 &  160 \\
        \hline
        $\dim (V)$ w/ RMCLI reduced, $q=3$ & 149 &  220 &  243  & 252 \\
        \hline
        $\dim (V)$ w/ RMCLI reduced, $q=5$ & 177  & 241 &  265  & 275 \\
        \hline
    \end{tabular}
    \caption{Dimensions of Rit subspace under RMCLI and its reduced version with $tol=10^{-6}$, $A_1,\ldots,A_d\in \mathbb{S}^{3969\times 3969}$ and $A_0,M\in \mathbb{S}^{3969\times 3969}_{++}$ from example \ref{example1}, Chebyshev degree $q=1,3,5$.}
    \label{tab:multi_dim_q}
\end{table}
Dimension of Ritz subspaces obtained from RMCLI and its reduced version for different $d$ and $q$ are collected in Table \ref{tab:multi_dim_q}. It is clear that both algorithms provide Ritz subspaces with small dimension and yield highly accurate Ritz eigenvalues. As a conclusion, we advise one to use RMCLI reduced with a tolerance of $10^{-5}$ or $10^{-6}$, small $q$ and $\varepsilon\approx \eta_1/10$ to reach sufficient precision in the cases studied here. 
\end{example}


\section{Conclusion}


In this paper, we provided a detailed study of a Ritz method for solving a class of multiparametric generalized eigenproblems based on sparse interpolation. We have shown that eigensolutions to \eqref{eq:evp} can be approximated using average matrices and a correction operator to build a good approximation Ritz subspace. Numerical experiments in section \ref{num_exp_section} show that the proposed method allows to approximate eigensolutions with high accuracy when the proposed Chebyshev-Legendre interpolation points are used. Based on results presented in section \ref{num_exp_section} a practical choice in theses examples is to compute the Ritz space is to use the reduced version of RMCLI with $q=3$ and $\varepsilon=\eta_1/10$ with $\eta=(\eta_1,\ldots,\eta_d)$ satisfying the condition to convergence of stochastic collocation interpolation. We did not address here the case of non-linear matrix operator $A:S\rightarrow \mathbb{S}^{n\times n}_{++}$ as the theory of average matrix could not be used so easily but could be done in the future to complete this study.

\bibliographystyle{siam}
\bibliography{biblio.bib}

\end{document}